%% file: Article.tex
\RequirePackage{fix-cm}
\documentclass[twocolumn]{svjour3}           % onecolumn (second format)
\smartqed  % flush right qed marks, e.g. at end of proof
\usepackage{graphicx}
\usepackage[switch,columnwise]{lineno} 	
\usepackage{natbib}
\usepackage{tabularx}
\usepackage{bm}
\usepackage{subcaption}
\usepackage{floatrow}
\usepackage{amsmath}
\usepackage{bm}
\usepackage{amssymb}
\usepackage{import}
\usepackage{booktabs}
\usepackage{dblfloatfix}
\usepackage{marvosym}
\usepackage[flushleft]{threeparttable}
\usepackage{listings}
\usepackage{color} 
\definecolor{mygreen}{RGB}{28,172,0} 
\definecolor{mylilas}{RGB}{170,55,241}
\usepackage{multirow}
\usepackage{verbatim}
\usepackage{rotating}
\usepackage{ulem}

\newcommand{\rfil}{r_\mathrm{fil}}
\newcommand{\rminsolid}[1]{r_\mathrm{min.Solid}^\mathrm{#1}}
\newcommand{\rminvoid}[1]{r_\mathrm{min.Void}^\mathrm{#1}}
\newcommand{\ProjField}[1]{{\bm{\bar{\rho}}}_\mathrm{#1}}
\newcommand{\Threshold}[1]{\eta_\mathrm{#1}}

\usepackage{xcolor} 
\definecolor{mygreen}{rgb}{0.0, 0.42, 0.24}

\journalname{PREPRINT}

\begin{document}
\title{Note on the minimum length scale and its defining parameters}
\subtitle{Analytical relationships for Topology Optimization based on uniform manufacturing uncertainties.}

\author{Denis Trillet \and
        Pierre Duysinx \and
        Eduardo Fern\'{a}ndez
}


\institute{Authors \at
           Department of Aerospace and Mechanical Engineering, \\ 
           University of Li\`{e}ge, 4000 Li\`{e}ge, Belgium. \\
							\\
		   Denis Trillet \\              
           \email{dtrillet@uliege.be}
}

\date{Received: date / Accepted: date}

\maketitle
%=============================================================================================
% |||||||||||||||||||||||||||||||||||||||||||||||||||||||||||||||||||||||||||||||||||||||||||
%=============================================================================================
\begin{abstract}
The robust topology optimization formulation that introduces the eroded and dilated versions of the design has gained increasing popularity in recent years, mainly because of its ability to produce designs satisfying a minimum length scale. Despite its success in various topology optimization fields, the robust formulation presents some drawbacks. This paper addresses one in particular, which concerns the imposition of the minimum length scale. In the density framework, the minimum size of the solid and void phases must be imposed implicitly through the parameters that define the density filter and the smoothed Heaviside projection. Finding these parameters can be time consuming and cumbersome, hindering a general code implementation of the robust formulation. Motivated by this issue, in this article we provide analytical expressions that explicitly relate the minimum length scale and the parameters that define it. The expressions are validated on a density-based framework. To facilitate the reproduction of results, MATLAB codes are provided.

As a side finding, this paper shows that to obtain simultaneous control over the minimum size of the solid and void phases, it is necessary to involve the 3 fields (eroded, intermediate and dilated) in the topology optimization problem. Therefore, for the compliance minimization problem subject to a volume restriction, the intermediate and dilated designs can be excluded from the objective function, but the volume restriction has to be applied to the dilated design in order to involve all 3 designs in the formulation.

\keywords{Length Scale \and Robust Design \and SIMP}
 
\end{abstract}
%=============================================================================================
% |||||||||||||||||||||||||||||||||||||||||||||||||||||||||||||||||||||||||||||||||||||||||||
%=============================================================================================

\section{Introduction} \label{sec:1}

Since the seminal work of \cite{BensoeKikuchi1988}, topology optimization has experienced huge advances and nowadays is being massively adopted in the industry \citep{Pedersen2006, Zhou2011, Zhu2016}. Among the successful advancements, one can mention the famous density method under the SIMP interpolation scheme \citep{Bendsoe1989}. The well-known shortcomings of SIMP led to a succession of improvements seeking to avoid the mesh dependency, the checkerboard patterns and the presence of intermediate densities. To date, one of the most effective approaches to dealing with the ill-effects of SIMP is the robust design approach that brings the eroded and dilated versions of the design \citep{Sigmund2009}. This method considers manufacturing errors that may incur in a uniformly thinner or uniformly thicker component compared to the blueprint layout. \citet{Sigmund2009} proposed a robust formulation that maximizes the performance of the worst performing design among the eroded, dilated and reference (intermediate) designs. This means, the formulation guarantees a design with good performance even if it is eventually eroded or dilated during the manufacturing process. Interestingly, the robust formulation yields a reference (intermediate) design that features minimum member size and minimum cavity size \citep{Wang2011}.

In addition to imposing minimum length scale in topology optimization, it has been seen that the robust formulation provides a more stable convergence than other projection or filtering strategies known so far, which allows to reach almost discrete solutions in the density method \citep{Wang2011}. The formulation has been recently applied in combination with existing methods that allow to impose stress limits \citep{DASILVA2019_2, daSilva2020}, overhang angle constraints \citep{Pellens2018}, maximum size restrictions \citep{Fernandez2020}, and geometric non-linearities \citep{Lazarov2011, daSilva2020}, not only in the density method but also in the level set method \citep{Chen2011,Andreasen2020}.

Despite proving its effectiveness in several fields of application \citep{Wang2011b, Christiansen2015}, the robust design approach presents some drawbacks with regard to its implementation. For example, the method requires boundary treatments with respect to the density filter, since the filtering region can be split at the edges of the design domain affecting the imposed minimum size \citep{Clausen2017,Kumar2021}. In addition, due to the erosion and dilation distances with respect to the intermediate design, boundary conditions can be disconnected from the designs involved in the formulation \citep{Clausen2017}. Another difficulty posed by the method is that the minimum length scale must be implicitly imposed through the parameters defining the density filter and the Heaviside projection. This drawback has been addressed in the literature using numerical \citep{Wang2011} and analytical \citep{Qian2013} approaches. The numerical approach consists of applying the density filter and the Heaviside projection to a 1D design, measuring the resulting length scale, and repeating the process several times with different filter and projection parameters to subsequently construct a graph that relates the involved parameters. The analytical method consists of applying the density filter as a convolution integral in a continuous 1D design. The solution of the integration yields explicit relationships of the projection/filtering parameters with the obtained length scale \citep{Qian2013}. As the three field scheme is a scalar function, the relationships obtained from the 1D designs remain valid for 2D and 3D \citep{Wang2011}. 

% =============================================================
% ====================== OLD ONE ==============================
% The numerical and theoretical relationships reported by \citet{Wang2011} and \citet{Qian2013} are intended for a particular case where the minimum size of the void phase is set equal to that of the solid phase. This simplifies and facilitates the procedures to relate the desired minimum length scale to the projection and filtering parameters. However, recent advances in robust formulation have highlighted the need to build a more extensive combination of parameters, in addition to accessing other relevant information. For example, \citet{Fernandez2020} show that in order to incorporate maximum size restrictions into the robust formulation, it is necessary to know the dilation and erosion distances, which until then had not been required by other methods.

The numerical and analytical relationships reported by \citet{Wang2011} and \citet{Qian2013} are intended for the particular case where the minimum size of the void phase is set equal to that of the solid phase. This simplifies and facilitates the procedures to explicitly relate the desired minimum length scale to the filter and projection parameters, but in this way, the minimum size of the solid phase cannot be different from the minimum size of the void phase. Nonetheless, due to the rapid popularization of the robust formulation, it becomes necessary to find a method that allows to impose different minimum sizes for each phase, in addition to provide access to other geometric information. For example, the erosion and dilation distances with respect to the intermediate design have been required to impose a maximum member size \citep{Fernandez2020} and to obtain designs tailored to the size of a deposition nozzle \citep{Fernandez2021}. To date, the erosion and dilation distances have been reported for specific minimum length scales, which hinders widespread applicability of the above methods.

% ==================================================
% ============== OLD ONE ===========================
%This work aims to extend the relationship between the desired length scale and the parameters that define it. The {\color{blue}purpose} is to obtain in a fast and precise way the filter and projection parameters that impose the user-defined minimum size. The study is not limited to a specific relationship between minimum size of the solid and void phases, and it also provides the dilation and erosion distances. To this end, we extend the analytical equations presented by \citet{Qian2013} and compare them with the numerical method proposed by \citet{Wang2011}. This comparison allows us to introduce correcting measures into the analytical method in order to consider the rounding errors coming from the finite element discretization.

The aim of this work is to provide a method to obtain the filter and projection parameters that impose user-defined minimum length scales, where the size of the solid phase can be defined differently from that of the void phase. In addition, for the desired minimum length scales, the erosion and dilation distances are provided. To this end, we extend the applicability of the analytical method proposed by \citet{Qian2013} and the resulting relationships are validated using the numerical method proposed by \citet{Wang2011} and a set of 2D topology-optimized designs. To facilitate the replication of results and the application of the methods discussed in this paper, MATLAB codes are provided.

%===============================================
%============ OLD ONE ==========================
%The remainder of the document is developed as follows. Section \ref{sec:2} introduces the robust formulation and the implemented case studies. Section \ref{sec:3} presents the analytical method for relating the minimum length scale to the parameters involved in the filtering and projection schemes. Section \ref{sec:4} compares the numerical and the analytical methods and provides the corrective measures to take into account the rounding errors. Section \ref{sec:5} presents some test cases to validate the obtained relationships, while the Section \ref{sec:6} gathers the final conclusions of this work. The codes required for the replication of results are provided in Section \ref{sec:7}.

The remainder of the document is developed as follows. Section \ref{sec:2} introduces the robust topology optimization formulation and the case studies used to asses the analytical expressions that relate the minimum len-gth scale to the parameters that impose it. Section \ref{sec:3} presents the analytical method proposed by \citet{Qian2013} and describes the main contribution of our work, which is the extension of the analytical approach to allow choosing independent minimum sizes for the solid phase and the void phase. Section \ref{sec:4} validates the analytical expressions using the numerical method proposed by \citet{Wang2011}. Section \ref{sec:5} discusses the sources of error that are inherent to the analytical method. Section \ref{sec:6} assesses the analytical expressions on 2D topology-optimized designs. Section \ref{sec:7} gathers the final conclusions of this work while Section \ref{sec:8} provides the codes that allow for the replication of results.

\section{Problem Definition} \label{sec:2}

The analytical expressions that relate the minimum length scale to the parameters that define it are developed for the robust topology optimization formulation based on the eroded, intermediate and dilated designs. This work considers the density approach based on the SIMP interpolation scheme \citep{Bendsoe1989}, even though the proposed methodology can be applied to other topology optimization approaches with little efforts.
 
Like most works in the literature, the eroded, intermediate and dilated designs that constitute the robust formulation are built using a three-field scheme \citep{Sigmund2013}. The first field, denoted by $\bm{\rho}$, corresponds to the design variables. The second field, denoted by $\bm{\tilde{\rho}}$, is obtained by a weighted average of the design variables within a circle of radius $r_\mathrm{fil}$. The third field, denoted by $\bm{\bar{\rho}}$, is obtained by projecting the components of the filtered field towards 0 or 1. The filter and projection functions are identical to those provided by \citet{Wang2011}, however, these are reminded herein for the sake of clarity as they define the minimum length scale obtained in the optimized design.

The filter of design variables, or the density filter \citep{Bruns2001,Bourdin2001}, is defined as follows:
\begin{equation}
    \tilde{\rho}_i = \frac{\displaystyle\sum_{j = 1}^N \rho_j \mathrm{v}_j w(\mathbf{x}_i-\mathbf{x}_j)}{\displaystyle\sum_{j = 1}^N {\mathrm{v}}_j w(\mathbf{x}_i-\mathbf{x}_j)} \;\; ,
\label{Eq:HatFilter}
\end{equation}
\noindent where $\rho_j$ is the design variable associated to the element $j$ and $\tilde{\rho}_i$ is the filtered variable associated to the element $i$. $\mathrm{v}_j$ is the volume of the element $j$ and $w(\mathbf{x}_i-\mathbf{x}_j)$ is the weigh of $\rho_j$ in the definition of $\tilde{\rho}_i$. As is common practice in the literature, the weighting function $w(\mathbf{x}_i-\mathbf{x}_j)$ is defined as a linear and decreasing function with respect to the distance between the elements $i$ and $j$, as follows:
\begin{equation}
    w(\mathbf{x}_i-\mathbf{x}_j) = \mathrm{max}\left(0,\: 1-\frac{\|\mathbf{x}_i-\mathbf{x}_j\|}{r_\mathrm{fil}}\right) \;\; ,
\label{Eq:Weighting_funtion}
\end{equation}
\noindent where $\mathbf{x}_i$ and $\mathbf{x}_j$ represent the centroid of the elements $i$ and $j$, respectively. It is recalled that $r_\mathrm{fil}$ is the radius of the density filter.

To reduce the amount of intermediate densities pre-sent in the filtered field and to build the designs that constitute the robust formulation, the projected field is obtained with the following smoothed Heaviside function \citep{Wang2011}:
\begin{equation}
\bar{\rho}_i = \frac{\tanh{(\beta \eta)} \tanh{(\beta(\Tilde{\rho}_i -\eta))}}
{\tanh{(\beta \eta)} \tanh{(\beta(1 - \eta))}}
\label{Eq:Heaviside}
\end{equation}
\noindent  where $\beta$ and $\eta$ control the steepness and the threshold of the projection, respectively. The eroded, intermediate and dilated designs, denoted by $\ProjField{ero}$, $\ProjField{int}$ and $\ProjField{dil}$, are obtained from the smoothed Heaviside function of Eq.~\eqref{Eq:Heaviside} by using the same $\beta$ but with different thresholds ${\Threshold{ero}}$, $\Threshold{int}$ and $\Threshold{dil}$, thus leading to ${\ProjField{ero}}(\bm{\tilde{\rho}},\beta,{\Threshold{ero}})$, ${\ProjField{int}}(\bm{\tilde{\rho}},\beta,{\Threshold{int}})$ and ${\ProjField{dil}}(\bm{\tilde{\rho}},\beta,{\Threshold{dil}})$.

\begin{figure}[t!]
	\captionsetup{width=1.00\linewidth}
    \centering	
    \begin{subfigure}{1.00\linewidth}
    	\centering
		\includegraphics[width=0.45\linewidth]{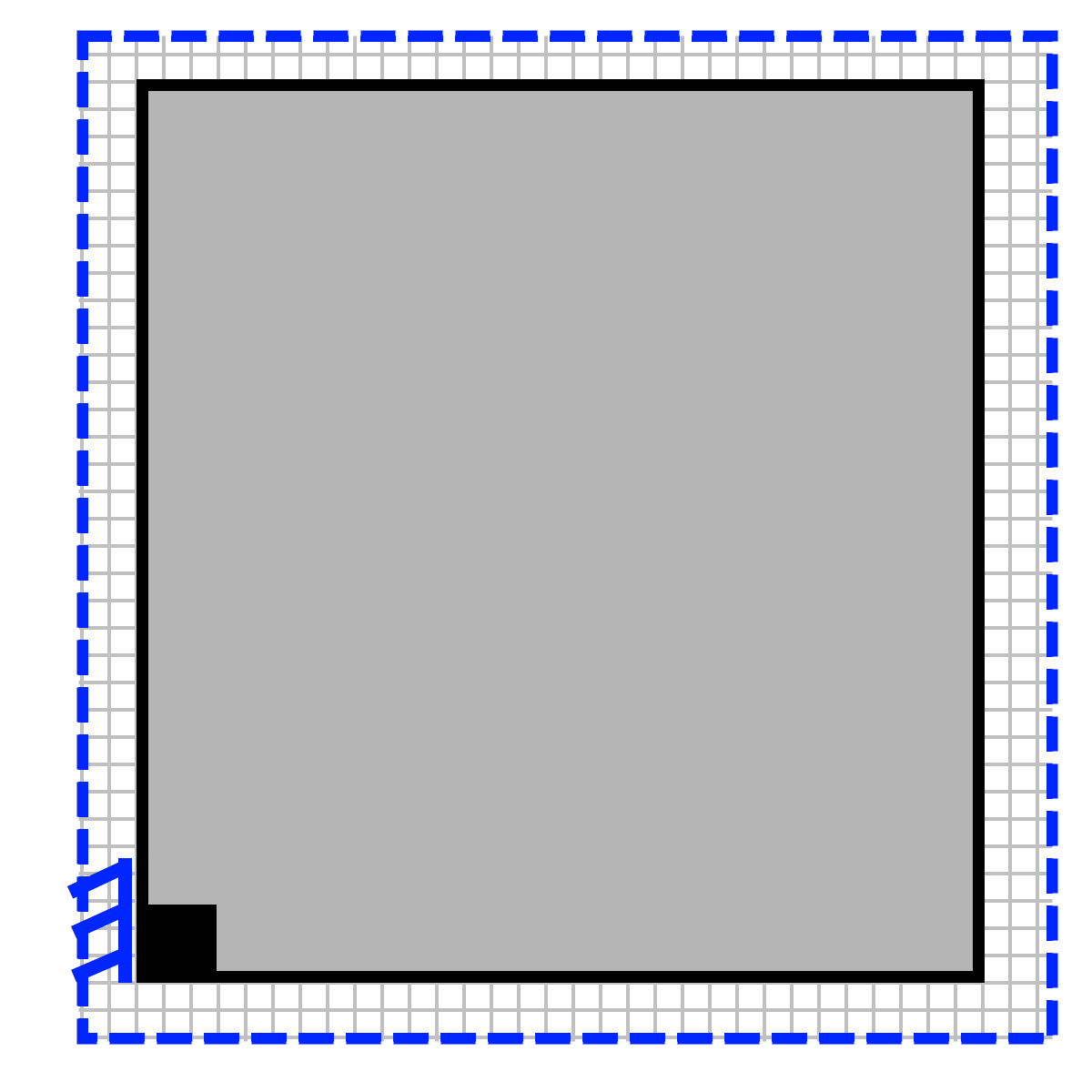}
	\end{subfigure}
	\vspace{2mm}\\
	\begin{subfigure}{1.00\linewidth}
    	\centering
		\includegraphics[width=0.70\linewidth]{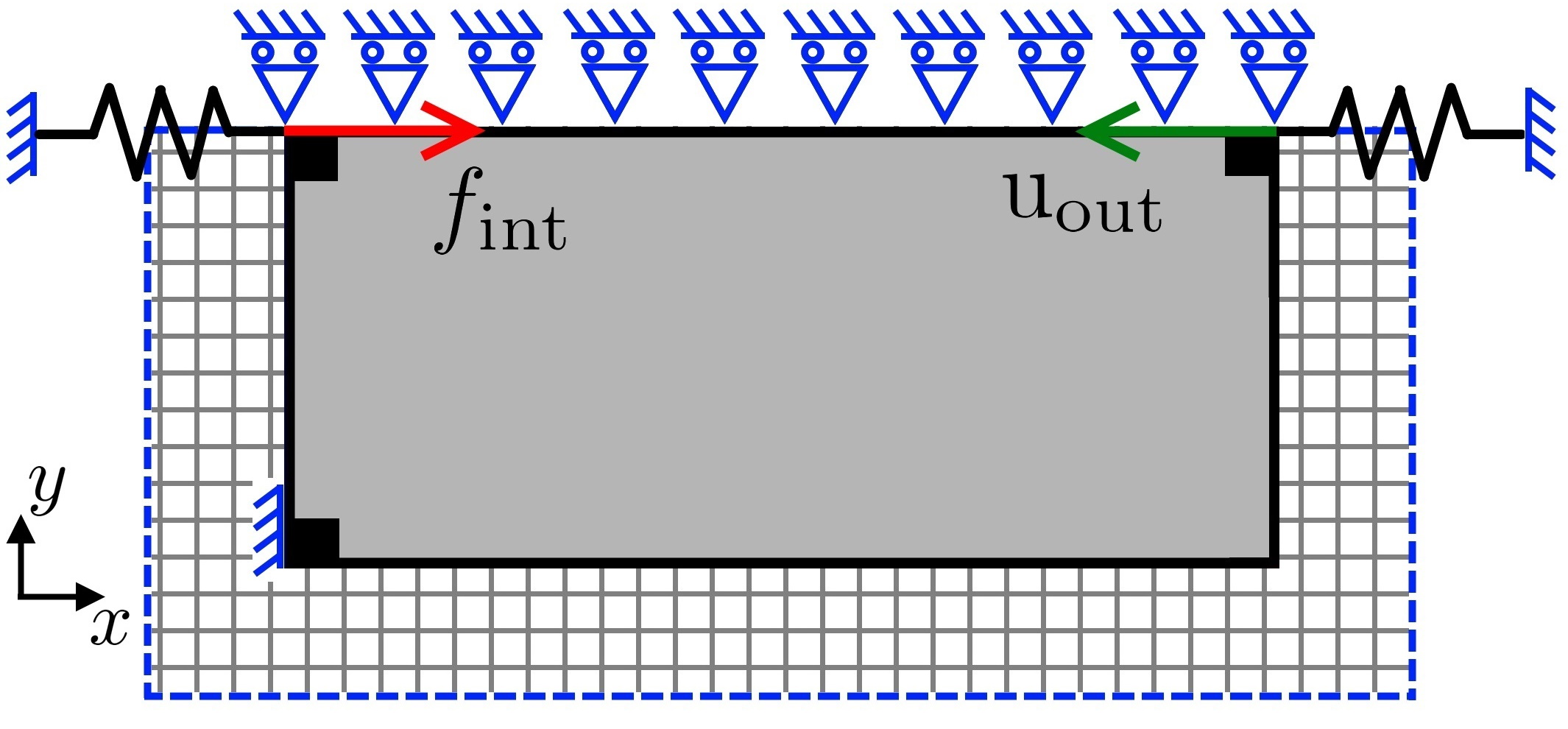}
	\end{subfigure}
	\vspace{-60mm}\\
	\hspace{-75mm} \footnotesize{(a)}
	\vspace{ 28mm}\\
	\hspace{-75mm} \footnotesize{(b)}
	\vspace{ 26mm}\\
	\caption{(a) Heat sink and (b) force inverter design domains considered in this study.}
	\label{fig:DesignDomains}
\end{figure}

To assess the scope of the equations that provide the necessary parameters to impose the minimum length scale, different topology optimization problems are sol-ved with variations in the desired minimum length scale. The minimum size of the optimized designs is measured and compared against the intended values. In this work, the problems that are chosen to illustrate the developments are the heat conduction and the non-linear compliant mechanism design. The former aims at minimizing the thermal compliance subject to a volume restriction \citep{Wang2011}, whilst in the latter, the formulation of a force inverter is considered where an output displacement is maximized for a given input force \citep{Sigmund1997}. The design domains are shown in Fig.~\ref{fig:DesignDomains}. According to the robust design approach, the topology optimization problems can be written as follows:

%From the optimized designs, 

%
\begin{equation}
\begin{split}
    \text{min}   \quad & \mathrm{max} \left( c({\ProjField{ero}}),\; c({\ProjField{int}})\;, c({\ProjField{dil}}) \right)   \\
    \text{s.t.:} \quad & \mathbf{v}^\intercal \: {\ProjField{dil}} \leq V^*_\mathrm{dil}(V^*_\mathrm{int}) \\
                & 0 \leq \rho_i \leq 1 \;\; , \;\; i=1,...,N\;\;,
\end{split}
\label{Eq:compMini}
\end{equation}
\noindent where $c$ represents the thermal compliance of the heat sink or the output displacement of the force inverter, $V^*_\mathrm{dil}$ is the upper bound of the volume restriction and N is the total amount of design variables. As the design intended for manufacturing is the intermediate design, the upper bound of the volume constraint is scaled according to the user-defined limit $V^*_\mathrm{int}$. The optimization problem and the scaling process of the volume restriction are the same as described by \cite{Wang2011} for the heat conduction problem. In the nonlinear force inverter, the method to deal with the mesh distorsion appearing in low density elements is the same as explained by \cite{Wang2014}. To avoid overextending the contents of the manuscript, the interested reader is refered to the cited articles.

\section{Analytical minimum length scale} \label{sec:3}

Comprehensive numerical tests have shown that the robust formulation in Eq.~\eqref{Eq:compMini} yields intermediate designs that have the same topology as their eroded and dilated projections \citep{Wang2011}. In other words, the erosion process does not destroy the solid members of the optimized topology (intermediate), it just makes them thinner. This implies that a minimum size of the solid members in the intermediate design is imposed according to the erosion distance. Similarly, the dilation process does not close the cavities of the optimized topology (intermediate), it only makes them narrower. This implies that the minimum size of the void phase is imposed according to the dilation distance. 

Considering a two-dimensional design domain, the minimum size of the solid phase is defined by the radius $r_\mathrm{min.Solid}$ of the largest circle that can be circumscribed into the smallest solid member of the topology, while the minimum size of the void phase is defined by the radius $r_\mathrm{min.Void}$ of the largest circle that can be circumscribed into the smallest cavity of the topology. It has been shown that the minimum length scale, $r_\mathrm{min.Solid}$ and $r_\mathrm{min.Void}$, are defined by the filter radius $r_\mathrm{fil}$ and the projection thresholds $\Threshold{ero}$, $\Threshold{int}$ and $\Threshold{dil}$. Therefore, if specific values are to be imposed for the minimum length scale, these must be implicitly imposed through the 4 parameters that define the density filter and the Heaviside projection. 

The dependence of the minimum length scale on the projection and filtering parameters hurdles a general implementation of the robust formulation, specially in commercial codes, as there is no prompt approach to find the relationship between the involved parameters and the desired length scale. Encouraged by this shortcoming, in this section we present and extend the applicability of the analytical procedure proposed by \citet{Qian2013}, which is described below.

\begin{figure}[t!]
	\captionsetup{width=1.00\linewidth}
    \centering	
    \begin{subfigure}{0.46\linewidth}
    	\centering
		\includegraphics[width=1.00\linewidth]{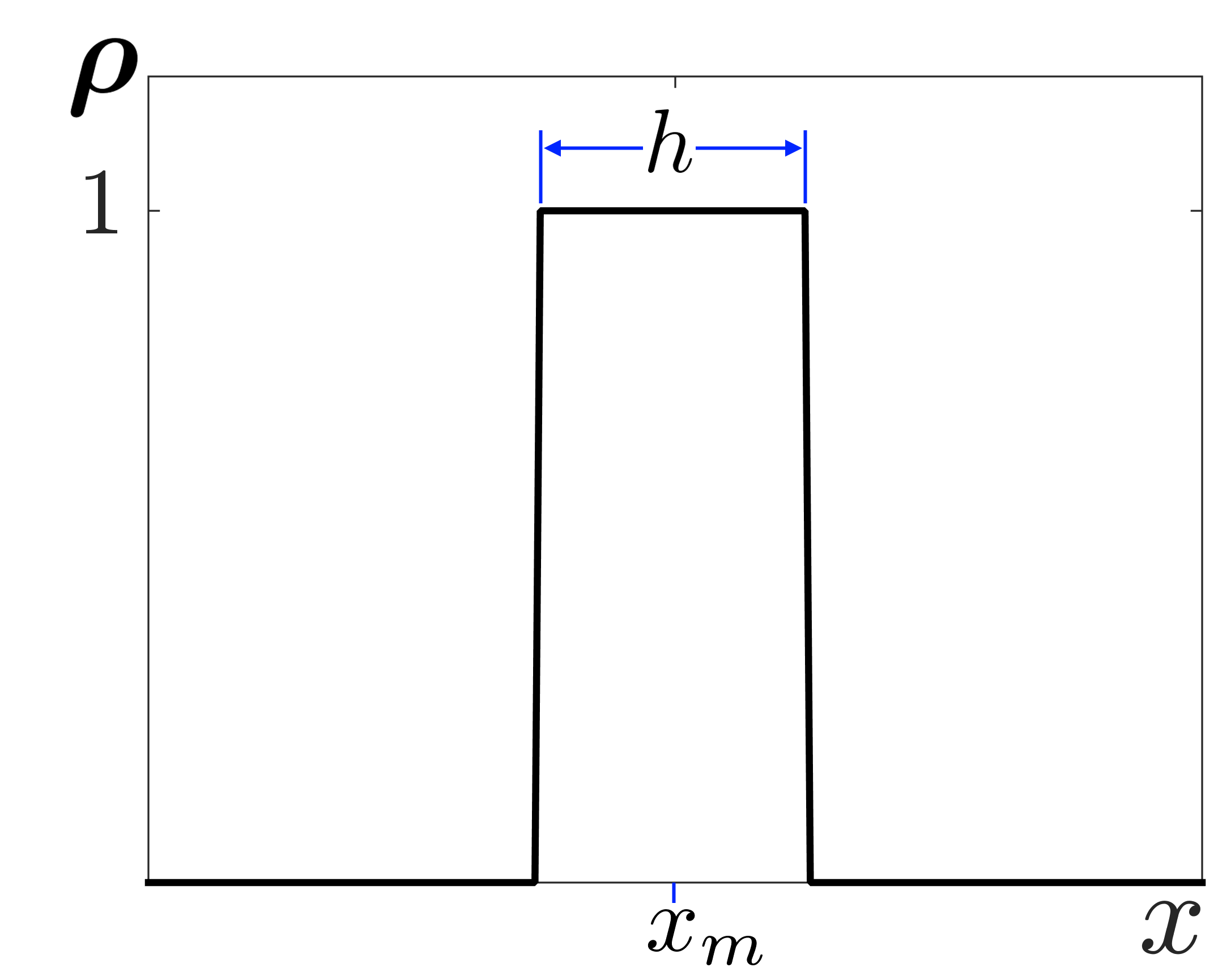}
		\caption{1D domain}
		\label{fig:fig_1_v2_a}
	\end{subfigure}
	~
	\begin{subfigure}{0.46\linewidth}
    	\centering
		\includegraphics[width=1.00\linewidth]{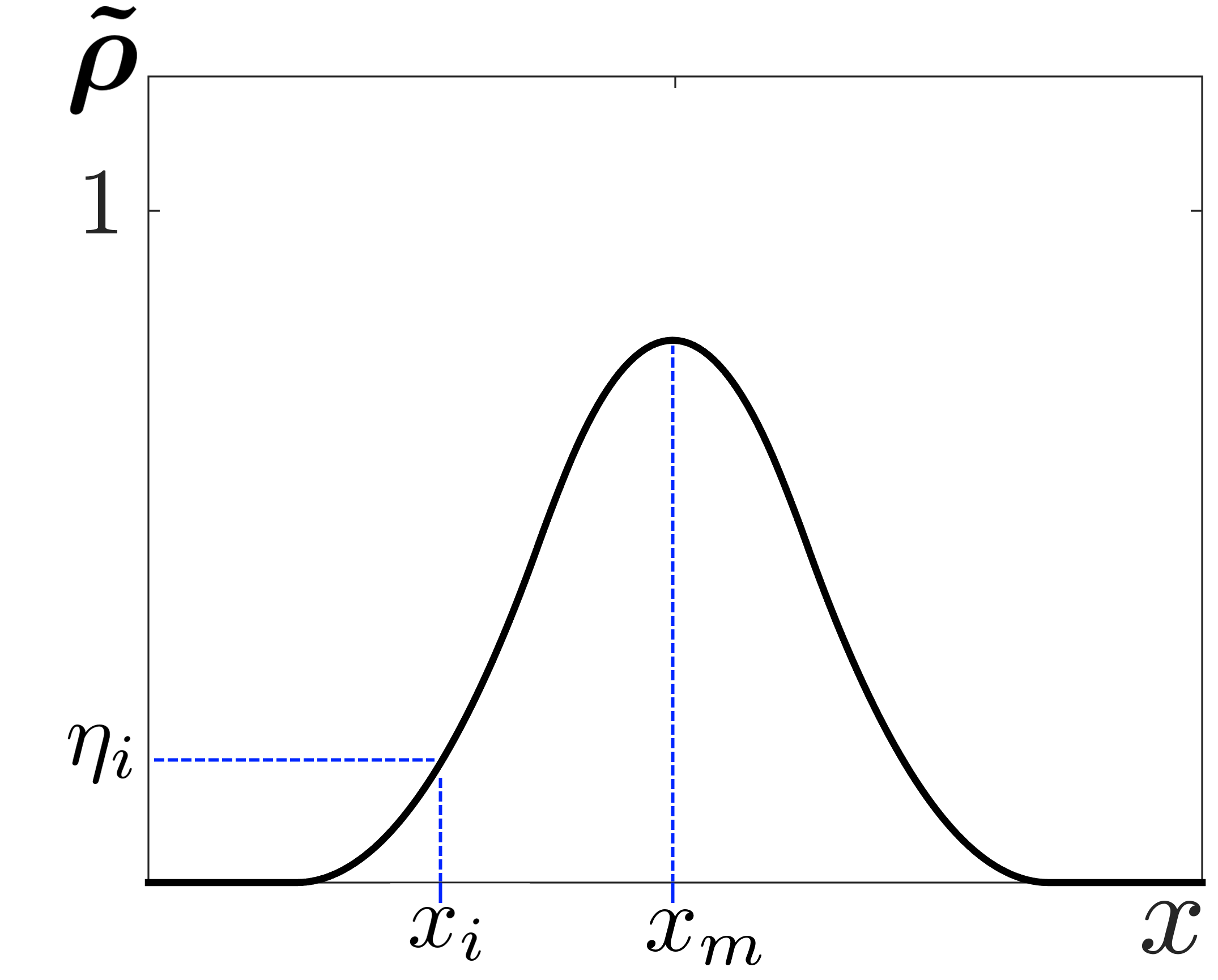}
		\caption{Filtered Field}
		\label{fig:fig_1_v2_b}
	\end{subfigure}
	\vspace{7mm}\\
	\begin{subfigure}{0.47\linewidth}
    	\centering
		\includegraphics[width=1.00\linewidth]{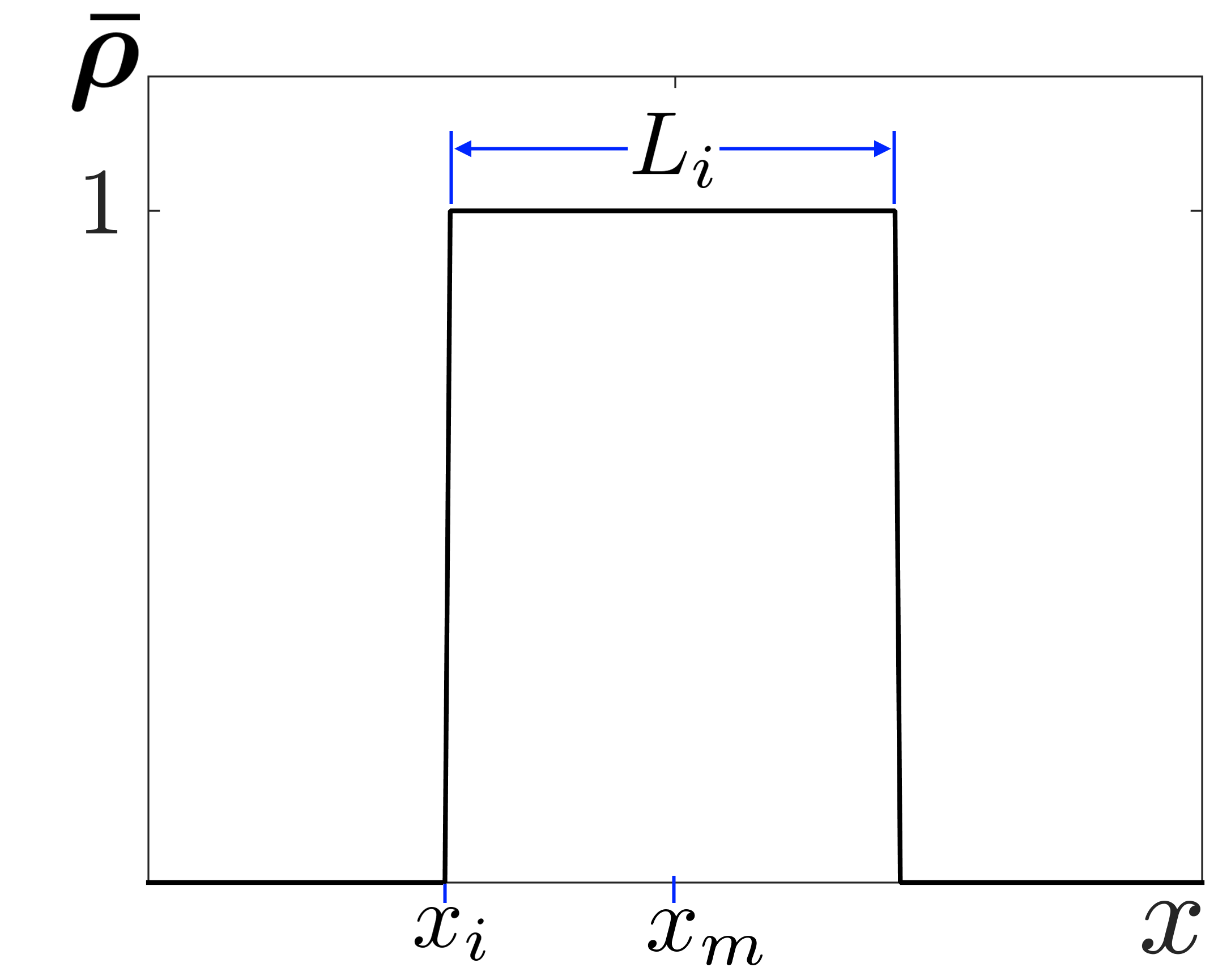}
		\caption{Projected Field}
		\label{fig:fig_1_v2_c}
	\end{subfigure}
	~
	\begin{subfigure}{0.47\linewidth}
    	\centering
		\includegraphics[width=1.00\linewidth]{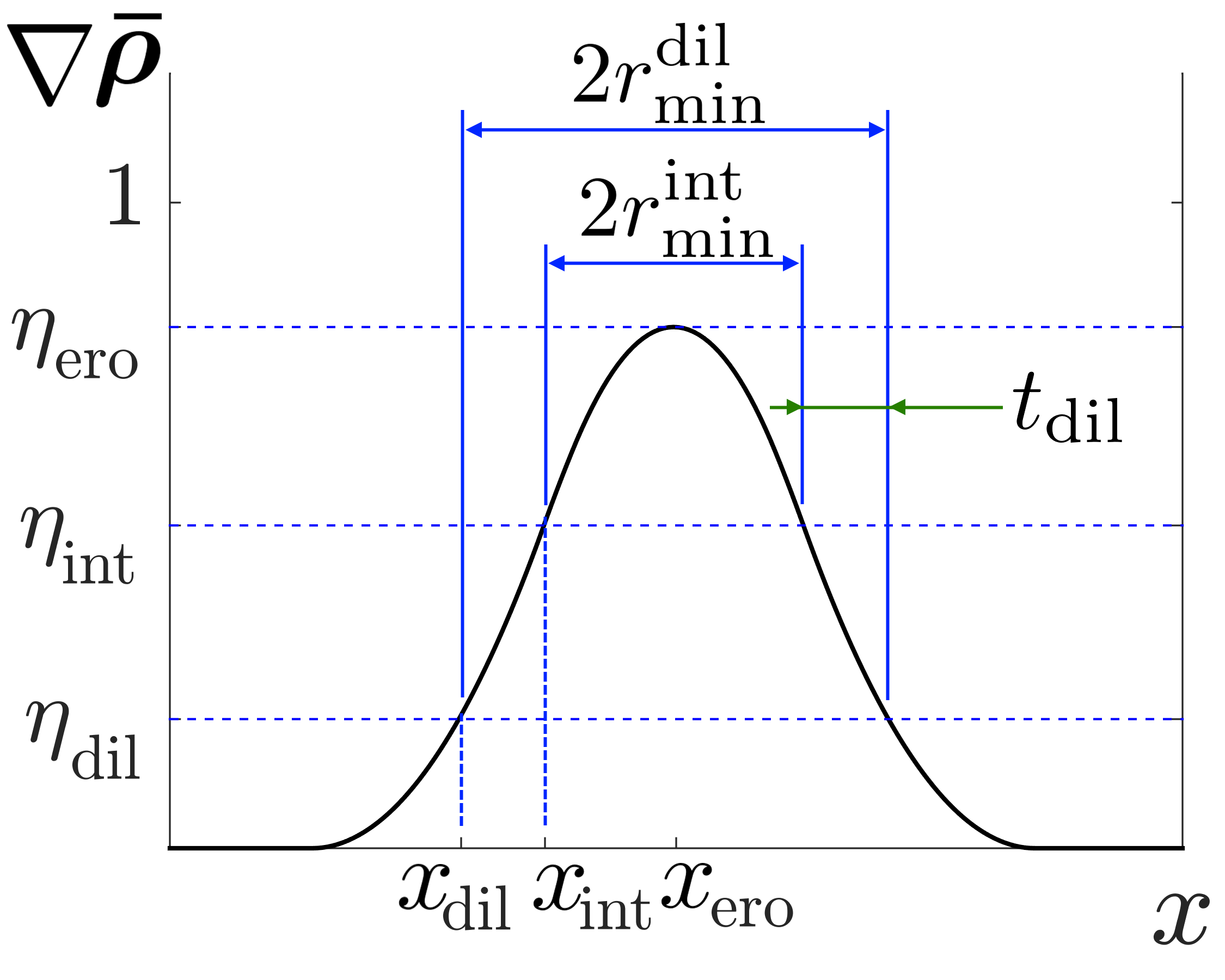}
		\caption{{Filtered field}}
		\label{fig:fig_1_v2_d}
	\end{subfigure}
	\caption{The three field scheme applied to a one-dimensional design domain in order to obtain the analytical minimum length scale.}
	\label{fig:fig_1_v2}
\end{figure}

To obtain an explicit relation between the minimum length scale and the parameters that define it, \citet{Qian2013} proposed to apply the three-field scheme over a uni-dimensional and continuous design domain. For instance, to obtain the minimum size of the solid phase, a design domain $\bm{\rho}(x)$ centered at the coordinate $x_m$ and containing solid elements on a stretch of size $h$ is assumed, as shown in Fig.~\ref{fig:fig_1_v2_a}. The continuous form of the density filter reads as follows:
\begin{equation}
\begin{matrix}
{\tilde{\rho}}(x_i) &=  \frac{\displaystyle\int_{x_i-r_\mathrm{fil}}^{x_i+r_\mathrm{fil}} \bm{\rho}(x) \left( 1-\frac{\left |x_i - x\right |}{r_\mathrm{fil}}\right) dx}
{\displaystyle\int_{x_i-r_\mathrm{fil}}^{x_i+r_\mathrm{fil}} \left( 1-\frac{\left |x_i - x \right |}{r_\mathrm{fil}}\right) dx }
\vspace{4mm}\\
&= \displaystyle\int_{x_i-r_\mathrm{fil}}^{x_i+r_\mathrm{fil}} \frac{\bm{\rho}(x)}{r_\mathrm{fil}} \left( 1-\frac{\left |x_i - x \right |}{r_\mathrm{fil}}\right) dx
\end{matrix}
\end{equation}
% xi is any coordinate?
\noindent where $x_i$ is any coordinate of $x$. For example, by choosing a filter radius greater than $h/2$, the one-dimensional filtered field in Fig.~\ref{fig:fig_1_v2_b} is obtained. The filtered field can be projected using the smoothed Heaviside function of Eq.~\eqref{Eq:Heaviside}. To simplify the analysis, an infinite steepness parameter ($\beta \to \infty$) is assumed. For example, for an Heaviside threshold $\eta_i=0.2$, the design of Fig.~\ref{fig:fig_1_v2_c} is obtained. The size of the solid phase in the projected design is defined by the length $L_i$, which can be obtained by finding $x_i$ from the equation $\tilde{\rho}(x_i)=\eta_i$. However, to this extent, the size of the solid phase is a function $L_i(r_\mathrm{fil},\eta_i,h)$ that depends on the assumed $h$ value. To determine the value of $h$ and to discover the explicit relation between the size $L$ of the projected field and the filter and projection parameters, it is necessary to refer to the foundation of the robust formulation.

An erosion projection removes solid material from the surface of the reference design. This operation enlarges the cavities and thins the structural members present in the intermediate reference design  \citep{Sigmund2009}. When robustness with respect to erosion is desired, the structural members must be present in both the eroded and the reference designs. It has been shown that a sufficient condition to preserve robustness is that the solid member in the eroded design has to be projected by at least an infinitesimal size, i.e.~$L({\Threshold{ero}})\approx0$, as shown in Fig~\ref{fig:fig_1_v2_d}. This condition allows determining the assumed size $h$ by solving the following equation:
\begin{equation}
\Tilde{\rho}(x_\mathrm{ero}) = \int_{-h/2}^{h/2} \frac{1}{r_\mathrm{fil}} \left(1-\frac{|x|}{r_\mathrm{fil}} \right) dx ={\Threshold{ero}}
\label{Eq:tildXp}
\end{equation}
\noindent where the reference system is placed at $x_\mathrm{ero}$ for integration. By solving the expression in Eq.~\eqref{Eq:tildXp}, the size $h$ is obtained as:
\begin{equation}
    h = 2 r_\mathrm{fil} (1 - \sqrt{1-{\Threshold{ero}}})
    \label{Eq:hSize}
\end{equation}

By obtaining the distance $h$ that produces an eroded projection of infinitesimal size, it is possible to relate the minimum size to the filter radius for any projection whose threshold meets $\eta_i<{\Threshold{ero}}$. A particular case is to choose $\eta_i=0.5$, which is the value chosen by \citet{Qian2013} to define the intermediate design. However, to extend the analytical method we use an arbitrary $\eta_i$ such that ${\Threshold{i}} < {\Threshold{ero}}$.

As mentioned above, the size of the solid phase $L$ in the projected field can be obtained by equating $\tilde{\rho}(x_i)$ to $\eta_i$. Since the condition of robustness is imposed on $h$ (Eq.~\ref{Eq:hSize}), the length corresponds to the minimum size, i.e.~$L=2r_\mathrm{min.Solid}$. The expression that allows to relate the minimum size to the filter radius is as follows:
\begin{equation} \label{Eq:Integral_for_L}
{\tilde{\rho}}(x_i) = \displaystyle\int_{x_i-r_\mathrm{fil}}^{x_i+r_\mathrm{fil}} \frac{\bm{\rho}(x)}{r_\mathrm{fil}} \left( 1-\frac{\left |x_i - x \right |}{r_\mathrm{fil}}\right) dx = \eta_i
\end{equation}

\begin{figure}[t!]
	\captionsetup{width=1.00\linewidth}
    \centering	
    \begin{subfigure}{0.48\linewidth}
    	\centering
		\includegraphics[width=1.0\linewidth]{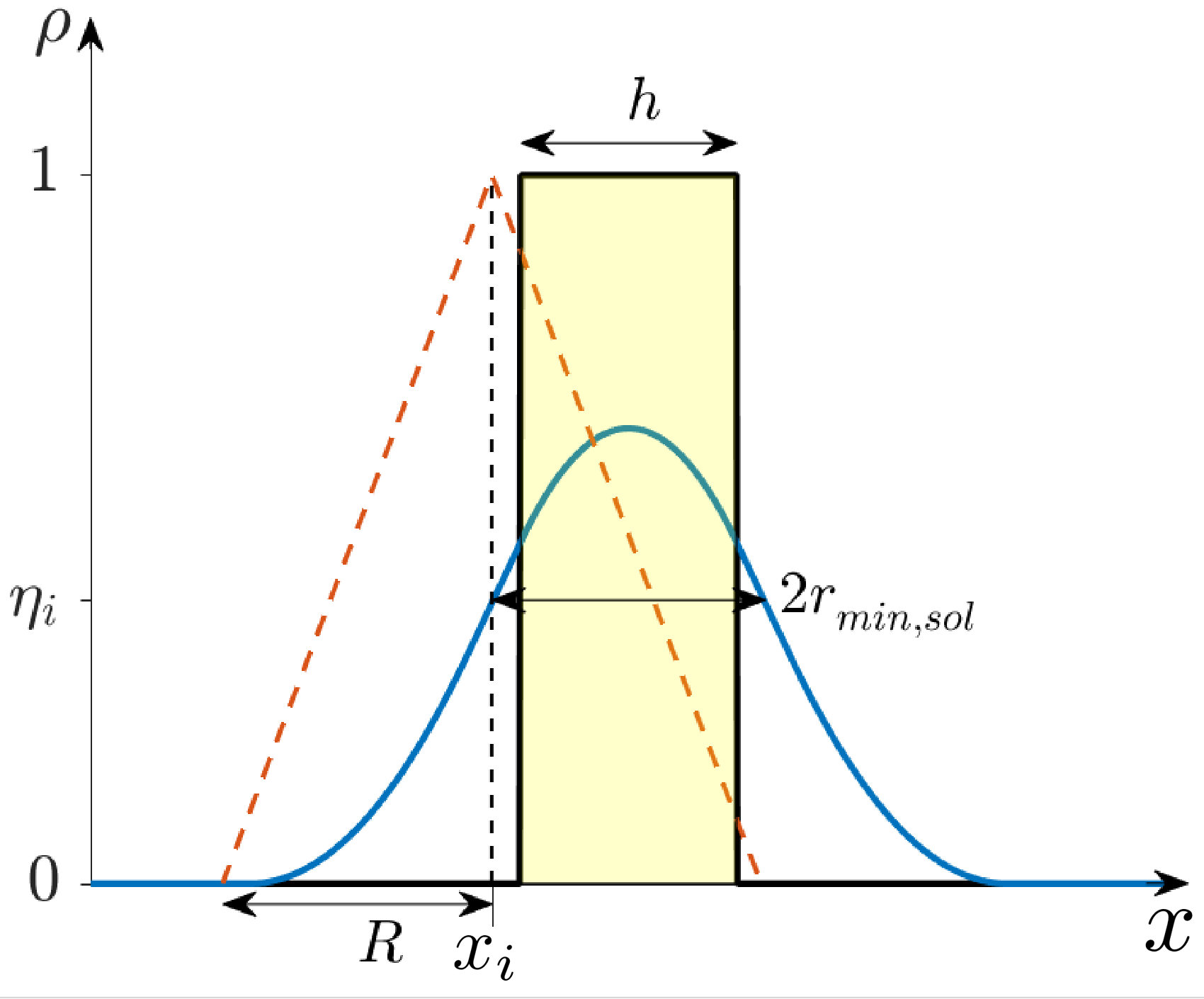}
		\caption{$L \geq h$ and $r_\mathrm{fil} \geq \frac{L+h}{2}$.}
		\label{fig:IntegralSize_v2_a}
	\end{subfigure}
    ~
	\begin{subfigure}{0.48\linewidth}
    	\centering
		\includegraphics[width=1.0\linewidth]{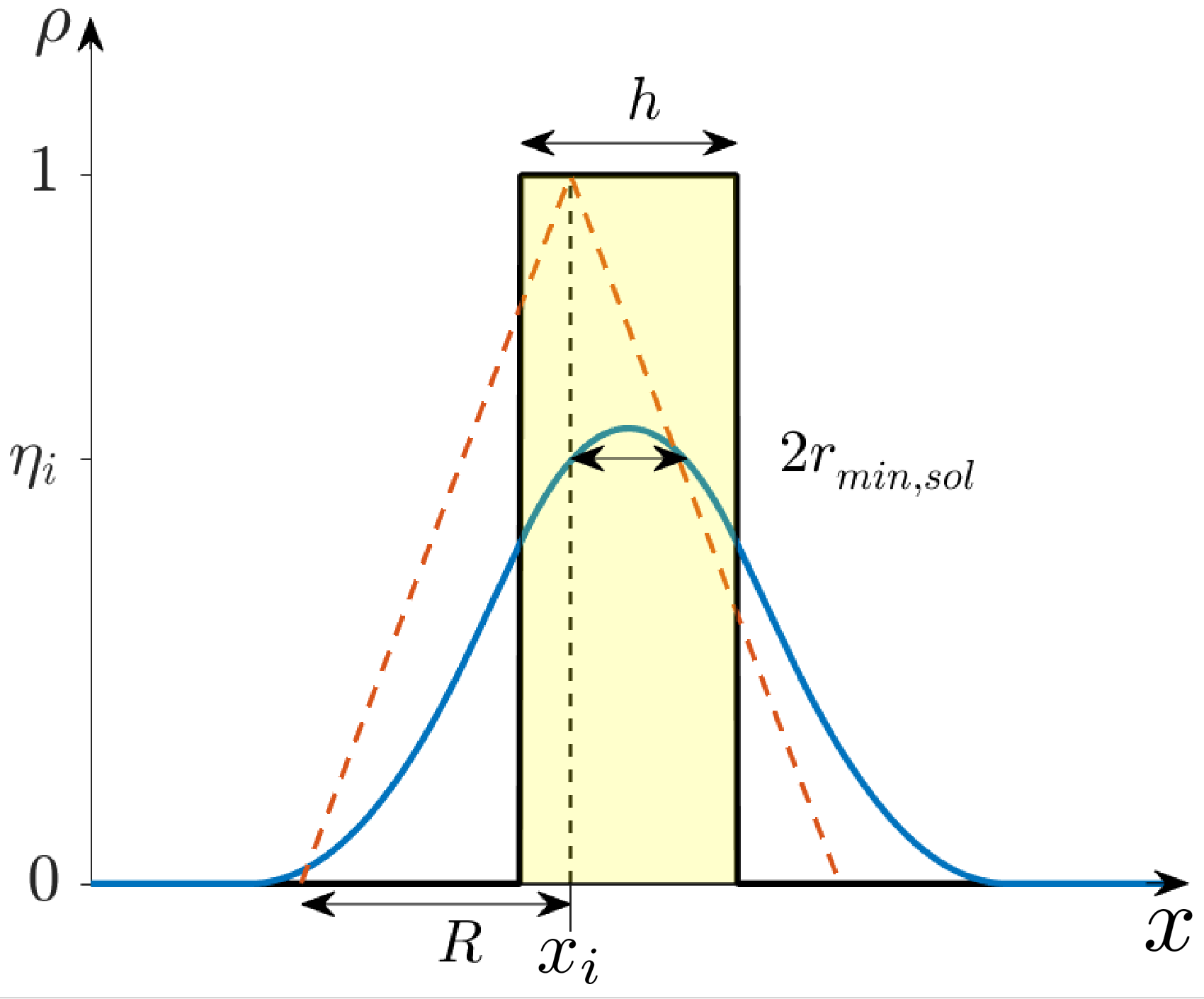}
		\caption{$L < h$ and $r_\mathrm{fil} \geq \frac{L+h}{2}$.}
		\label{fig:IntegralSize_v2_b}
	\end{subfigure}
	\vspace{3mm}\\
	\begin{subfigure}{0.47\linewidth}
    	\centering
		\includegraphics[width=1.0\linewidth]{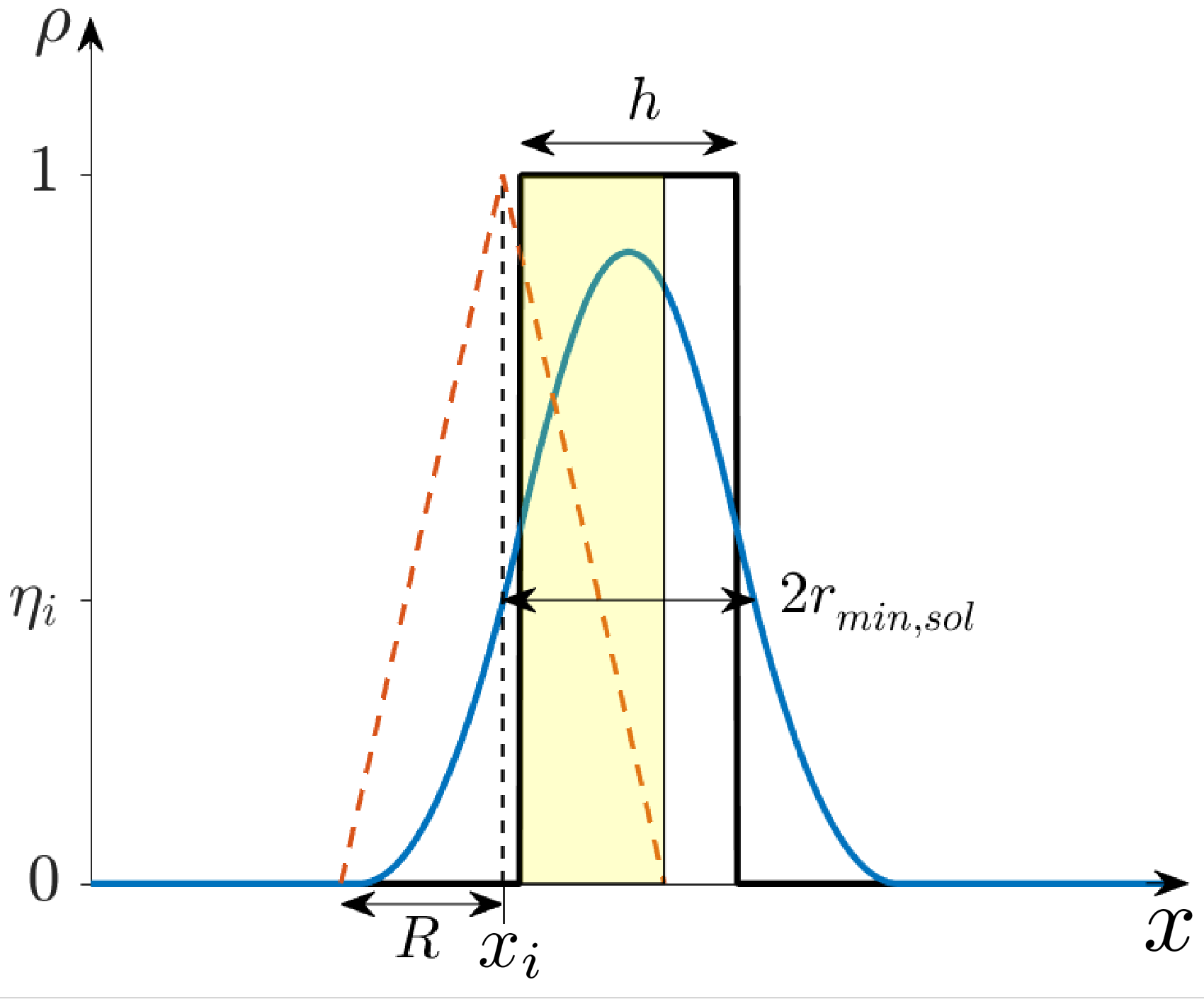}
		\caption{$L \geq h$ and $r_\mathrm{fil} < \frac{L+h}{2}$.}
		\label{fig:IntegralSize_v2_c}
	\end{subfigure}
	~
	\begin{subfigure}{0.45\linewidth}
    	\centering
		\includegraphics[width=1.0\linewidth]{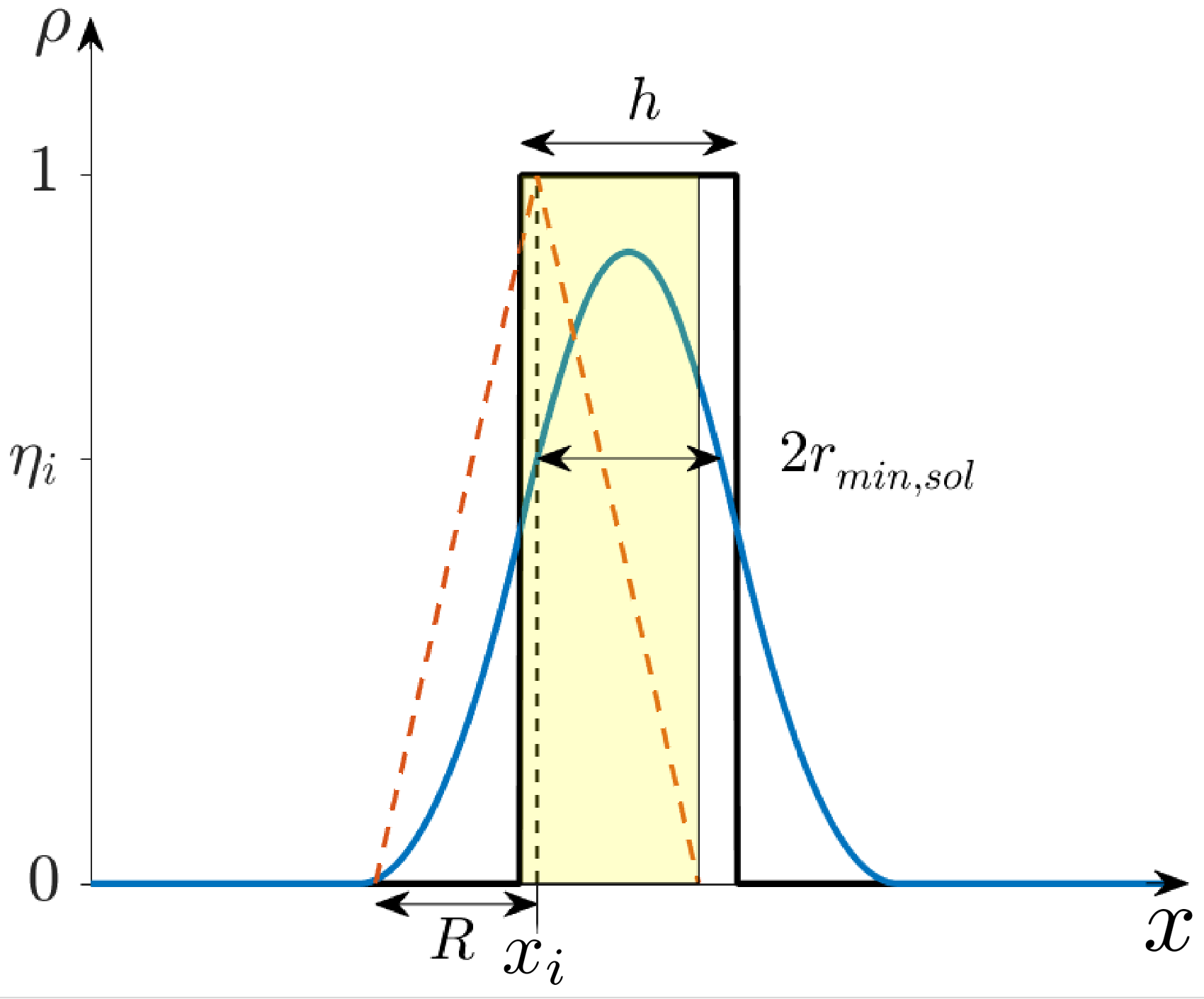}
		\caption{$L < h$ and $r_\mathrm{fil} < \frac{L+h}{2}$.}
		\label{fig:IntegralSize_v2_d}
	\end{subfigure}
	\caption{Situations to be considered when defining the limits of integration.}
	\label{fig:IntegralSize_v2}
\end{figure}

For the integration, it is convenient to place the origin of the reference system at $x_i$. To this end, 4 situations must be considered to define the integration limits, which are summarized in Fig.~\ref{fig:IntegralSize_v2}. For example, if the length of the projected field $L$ is greater than $h$ and the filter radius is greater than $(L+h)/2$ (Fig.~\ref{fig:IntegralSize_v2_a}), then Eq.~\eqref{Eq:Integral_for_L} becomes:
\begin{equation} \label{Eq:Integral_for_L_part_2}
{\tilde{\rho}}(x_i) = \displaystyle\int_{\frac{L-h}{2}}^{L-\frac{L-h}{2}} \frac{1}{r_\mathrm{fil}} \left( 1-\frac{\left | x \right |}{r_\mathrm{fil}}\right) dx = \eta_i
\end{equation}

%========================================================================
%================= FIGURE ===============================================
%========================================================================
%========================================================================

%========================================================================
%==================== TABLE =============================================
\begin{table*}
    \setlength\extrarowheight{4pt}
	\captionsetup{width=0.99\linewidth}
	\renewcommand{\arraystretch}{1.3}
    \centering
    \begin{tabular}{c r c l c l}
        \toprule
        Phase & \multicolumn{3}{c}{Conditions} & & \hspace{15mm}Minimum size 
        \\
        \hline 
		\multirow{4}{*}[-2em]{{\rotatebox{90}{\large{Solid}}}}
		&        
        $\eta_i \leq 2 {\Threshold{ero}} - 1$
        &
        and
        &
        $\eta_i \geq 0.5$ 
        & & 
        $\displaystyle\frac{2r_\mathrm{min.Solid}}{r_\mathrm{fil}} = 2 \sqrt{2 - 2 \eta_i} - 2 \sqrt{1 - {\Threshold{ero}}}$ 
        \vspace{1mm}\\  
        \cmidrule{2-6}  
        &  
        $\eta_i \geq 2 {\Threshold{ero}} - 2 + 2 \sqrt{1-{\Threshold{ero}}}$
		&        
        and
        &
        $\eta_i > 2 {\Threshold{ero}} - 1$
        & & 
        $\displaystyle\frac{2r_\mathrm{min.Solid}}{r_\mathrm{fil}} = 2 \sqrt{{\Threshold{ero}} - \eta_i}$ 
        \vspace{1mm}\\  
        \cmidrule{2-6}      
        &
        $\eta_i < 4 - 4 \sqrt{1-{\Threshold{ero}}} - 2{\Threshold{ero}}$
        &
        and
        & 
        $\eta_i < 0.5$ 
        & & 
        $\displaystyle\frac{2r_\mathrm{min.Solid}}{r_\mathrm{fil}} = 4 - 2 \sqrt{1 - {\Threshold{ero}}} - 2 \sqrt{2 \eta_i}$ 
        \vspace{1mm}\\  
        \cmidrule{2-6} 
        &   
        $\eta_i \geq 4 - 4 \sqrt{1-{\Threshold{ero}}} - 2{\Threshold{ero}}$
        &
        and
        &
        $\eta_i < 2 \Threshold{ero} - 2 + 2 \sqrt{1-{\Threshold{ero}}}$
        & & 
        $\displaystyle\frac{2r_\mathrm{min.Solid}}{r_\mathrm{fil}} = 2 - \frac{\eta_i}{1-\sqrt{1-{\Threshold{ero}}}}$ 
        \vspace{1mm}\\
        \hline 
		\multirow{4}{*}[-2em]{{\rotatebox{90}{\large{Void}}}}
		&        
        $\eta_i \geq 2 {\Threshold{dil}}$
        &
        and
        &  
        $\eta_i \leq 0.5$ 
        & &
        $\displaystyle\frac{2r_\mathrm{min.Void}}{r_\mathrm{fil}} = 2 \sqrt{2 \eta_i} - 2 \sqrt{{\Threshold{dil}}}$ 
        \vspace{1mm}\\  
        \cmidrule{2-6}  
        &  
        $\eta_i \leq 2 \eta_d + 1 - 2 \sqrt{{\Threshold{dil}}}$
		&        
        and
        &
        $\eta_i < 2 {\Threshold{dil}}$
        & & 
        $\displaystyle\frac{2r_\mathrm{min.Void}}{r_\mathrm{fil}} = 2 \sqrt{\eta_i - {\Threshold{dil}}}$ 
        
        \vspace{1mm}\\  
        \cmidrule{2-6}      
        &
        $\eta_i \geq 4 \sqrt{{\Threshold{dil}}} - 2{\Threshold{dil}} - 1$
        &
        and
        & 
        $\eta_i > 0.5$ 
        & & 
        $\displaystyle\frac{2r_\mathrm{min.Void}}{r_\mathrm{fil}} = 4 - 2 \sqrt{ {\Threshold{dil}}} - 2 \sqrt{2 - 2 \eta_i}$ 
        
        \vspace{1mm}\\  
        \cmidrule{2-6} 
        &   
        $\eta_i < 4 \sqrt{{\Threshold{dil}}} - 2{\Threshold{dil}} - 1$
        &
        and
        &
        $\eta_i > 2 {\Threshold{dil}} + 1 - 2 \sqrt{{\Threshold{dil}}}$
        & & 
        $\displaystyle\frac{2r_\mathrm{min.Void}}{r_\mathrm{fil}} = 2 - \frac{1-\eta_i}{1-\sqrt{{\Threshold{dil}}}}$
        \\  
        \bottomrule	
    \end{tabular}
    \caption{The explicit relationship between the minimum length scale and the filter and projection parameters.}
    \label{tab:SolidZoneCondition_v2}
\end{table*}
%========================================================================
%========================================================================

Solving the integral of Eq.~\eqref{Eq:Integral_for_L_part_2} leads to the following expression:

\begin{equation} \label{Eq:Integral_for_L_part_3}
{\tilde{\rho}}(x_i) = \frac{h}{r_\mathrm{fil}}(1-\frac{L}{2 r_\mathrm{fil}}) = \eta_i
\end{equation}

Finally, replacing Eq.~\eqref{Eq:hSize} in \eqref{Eq:Integral_for_L_part_3}:

\begin{equation} \label{Eq:Integral_for_L_part_4}
\frac{2r_\mathrm{min.Solid}}{r_\mathrm{fil}} = 2 -  \frac{\eta_i}{1-\sqrt{1-{\Threshold{ero}}}}
\end{equation}

Eq.~\eqref{Eq:Integral_for_L_part_4} explicitly relates the minimum size of the solid phase ($r_\mathrm{min.Solid}$), of a projected field defined by $\eta_i$, with the filter radius ($r_\mathrm{fil}$) and the erosion threshold ($\eta_\mathrm{ero}$). However, Eq.~\eqref{Eq:Integral_for_L_part_4} is only valid for $L> h$ and  $r_\mathrm{fil} \geq (L+h)/2$. For implementation purposes, it is more convenient to express the range of application in terms of the projection thresholds. To this end, $h$ can be replaced from Eq.~\eqref{Eq:hSize} and $L$ from \eqref{Eq:Integral_for_L_part_4}, which leads to conditions depending only on $\eta_i$ and ${\Threshold{ero}}$, as follows:
\begin{equation} \label{Eq:Integral_for_L_part_5}
\begin{matrix}
L > h & \implies \eta_i < 2 \Threshold{ero} - 2 + 2 \sqrt{1-{\Threshold{ero}}} \\[1ex]
r_\mathrm{fil} \geq (L+h)/2 &\implies \eta_i \geq 4 - 4 \sqrt{1-{\Threshold{ero}}} - 2{\Threshold{ero}}
\end{matrix}
\end{equation}

By repeating the procedure from Eq.~\eqref{Eq:Integral_for_L_part_2} to \eqref{Eq:Integral_for_L_part_5} for the 4 integration conditions shown in Fig.~\ref{fig:IntegralSize_v2}, a set of equations is obtained which relate the filter and projection parameters with the minimum size for any projection threshold $\eta_i$, provided that $\eta_i<{\Threshold{ero}}$. The set of equations are summarized in the four first rows of Table \ref{tab:SolidZoneCondition_v2}. 

To obtain the relationships that define the minimum size of the void phase, the same procedure must be used as for the solid phase, however, now starting from a one-dimensional design domain containing a cavity of size $h$. To avoid overextending the document with redundant information, this section is limited to presenting the final equations that define the minimum size of the void phase. The expressions are summarized in the last 4 rows of Table \ref{tab:SolidZoneCondition_v2}.

Having delivered the set of equations that expand the scope of the method proposed by \citet{Qian2013}, the following section presents a methodology to use these equations. 

\section{Imposing the desired minimum length scale}\label{sec:4}

In structural design, the minimum length scale control is usually desired because of design requirements or manufacturing limitations, hence in most cases, the minimum size of the solid and void phases are known values established for the intermediate design. Therefore, for the set of equations presented in Table \ref{tab:SolidZoneCondition_v2}, the radii $r_\mathrm{min.Solid}^\mathrm{int}$ and $r_\mathrm{min.Void}^\mathrm{int}$ are assumed user-defined input values. In this case, the projection threshold \linebreak $\eta_i$ corresponds to the projection threshold ${\Threshold{int}}$, \linebreak hence the desired length scale for the intermediate design is a function of the projection thresholds and of the size of the filter, namely, $r_\mathrm{min.Solid}^\mathrm{int}({\Threshold{int}}, {\Threshold{ero}}, r_\mathrm{fil})$ and $r_\mathrm{min.Void}^\mathrm{int}({\Threshold{int}}, {\Threshold{dil}}, r_\mathrm{fil})$. Given the number of unknowns $(r_\mathrm{fil}, {\Threshold{ero}}, \Threshold{int}, \Threshold{dil})$, the system of equations in Table \ref{tab:SolidZoneCondition_v2} becomes indeterminate and the desired minimum length scale can be imposed through multiple combinations of parameters. Nevertheless, such freedom of parameters selection can be reduced by considering the following three recommendations.

\begin{figure}[t!]
	\captionsetup{width=1.00\linewidth}
    	\centering
		\includegraphics[width=0.95\linewidth]{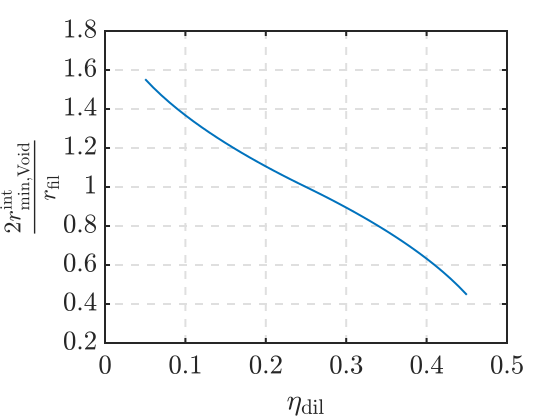}
	\caption{ Graphical relationship between the minimum size of the void phase, the filter radius and the dilation projection.}
	\label{fig:MinVoidFilter}
\end{figure}

\begin{figure*}[!p]
	\captionsetup[subfigure]{labelformat=empty}
	\captionsetup{width=1.00\linewidth}
    \centering
    \hspace{2mm}	
    \begin{subfigure}{0.48\linewidth}
    	\centering
		\includegraphics[width=0.92\linewidth]{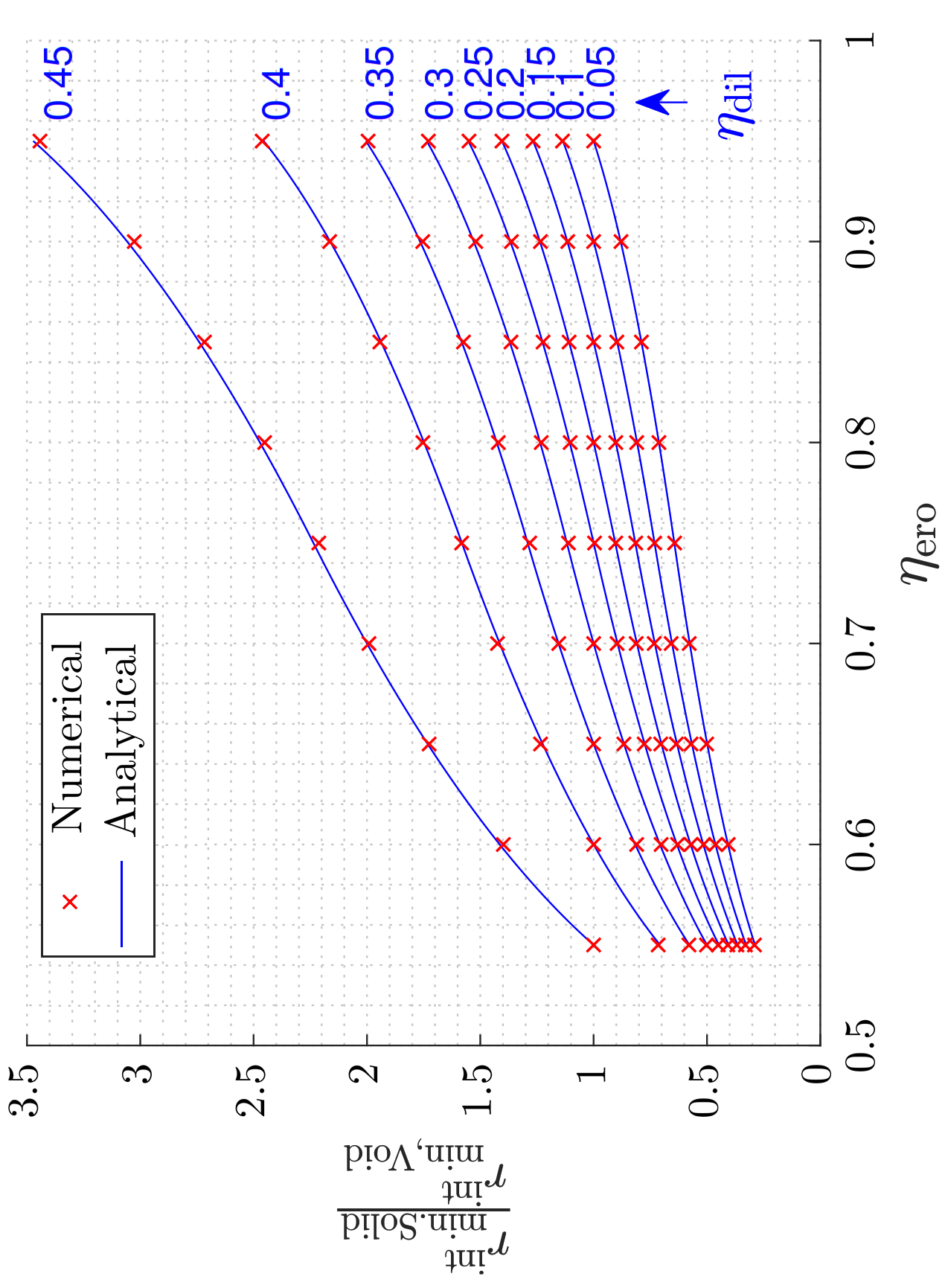}
		\caption{}
		\label{fig:fig_1_a}
	\end{subfigure}
	~
	\begin{subfigure}{0.48\linewidth}
    	\centering
		\includegraphics[width=0.92\linewidth]{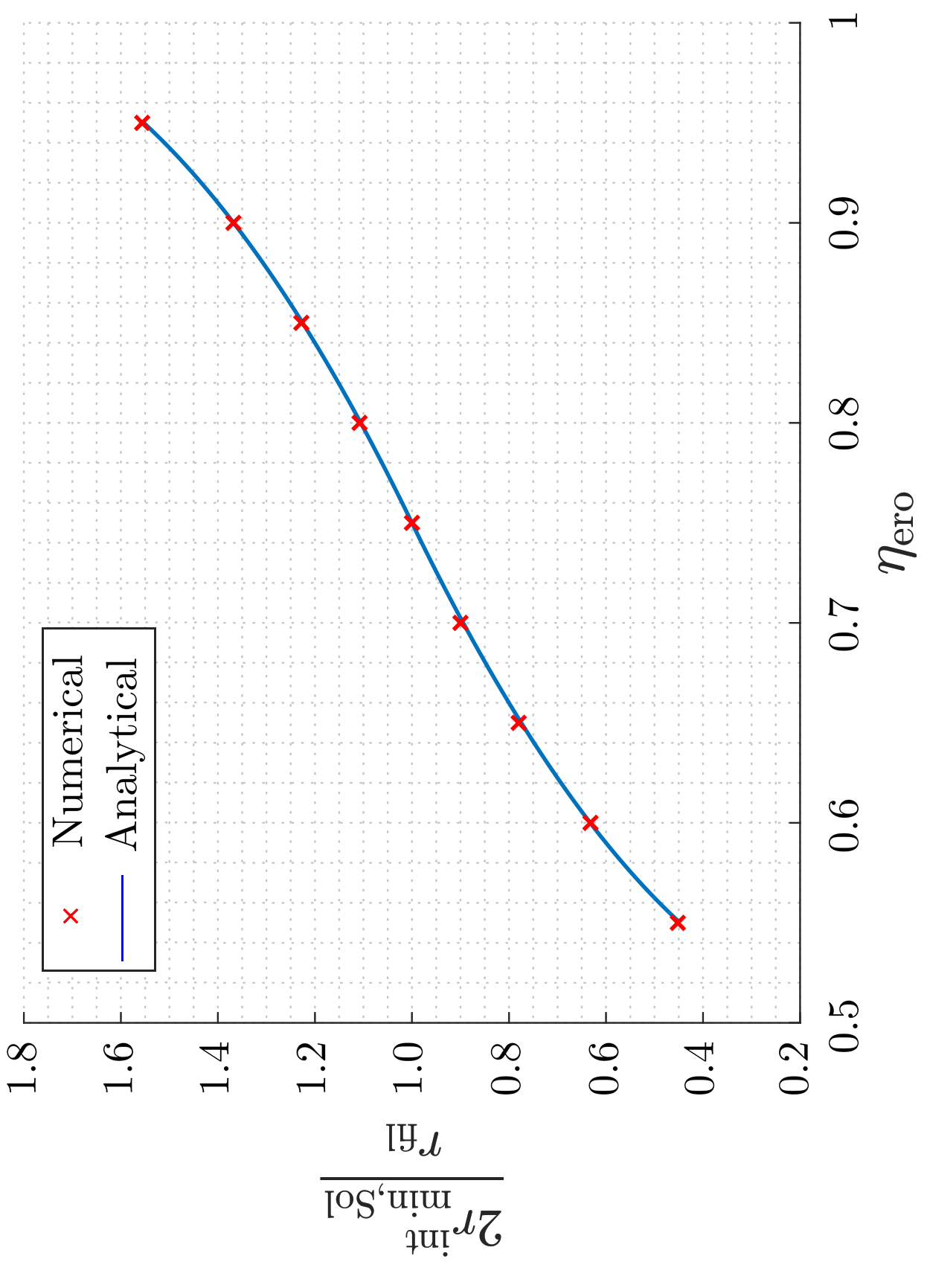}
		\caption{}
		\label{fig:fig_1_b}
	\end{subfigure}
	\vspace{-3mm}\\
	\hspace{2mm}	
	\begin{subfigure}{0.48\linewidth}
    	\centering
		\includegraphics[width=0.92\linewidth]{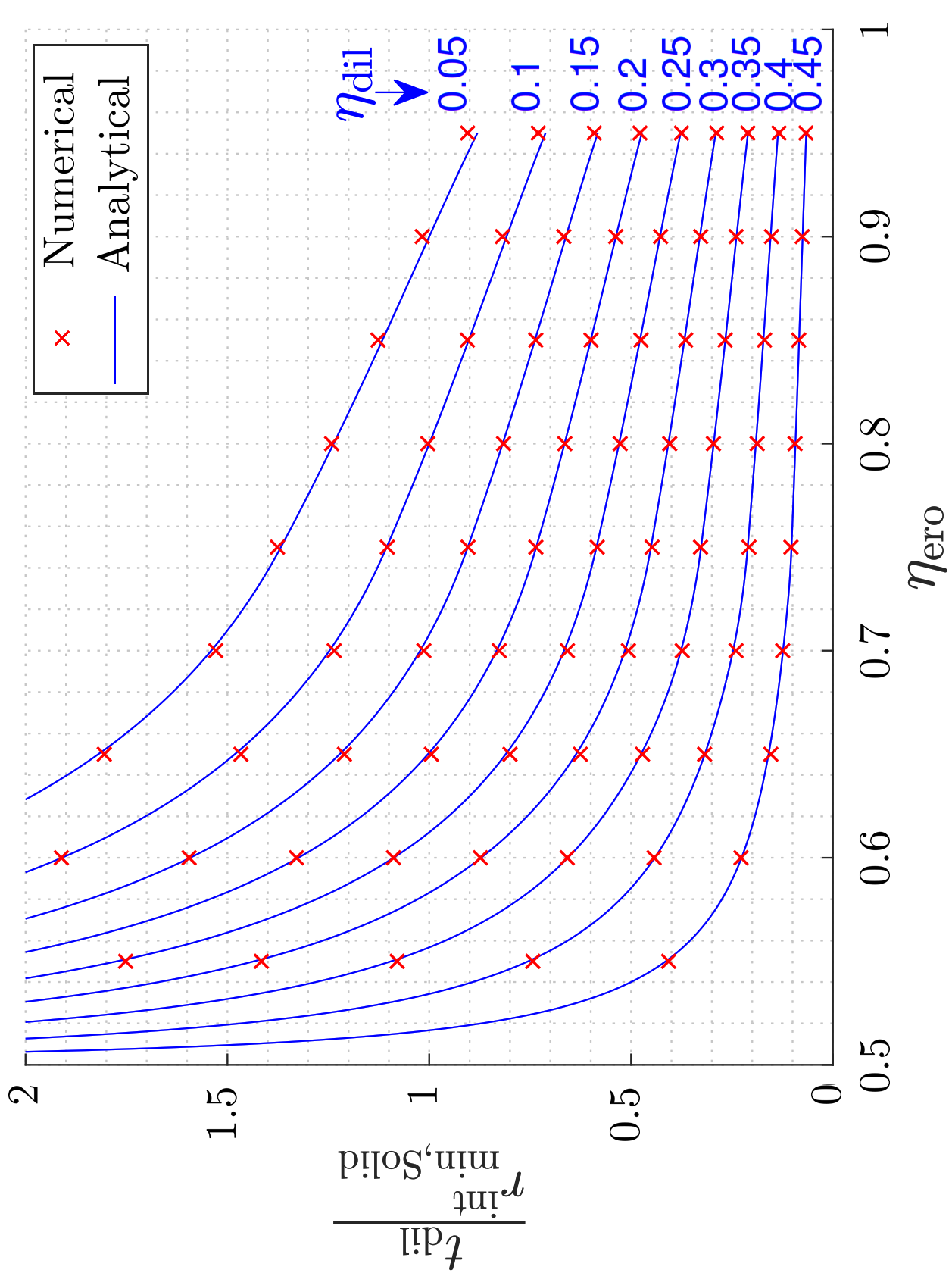}
		\caption{}
		\label{fig:fig_1_c}
	\end{subfigure}
	~
	\begin{subfigure}{0.48\linewidth}
    	\centering
		\includegraphics[width=0.92\linewidth]{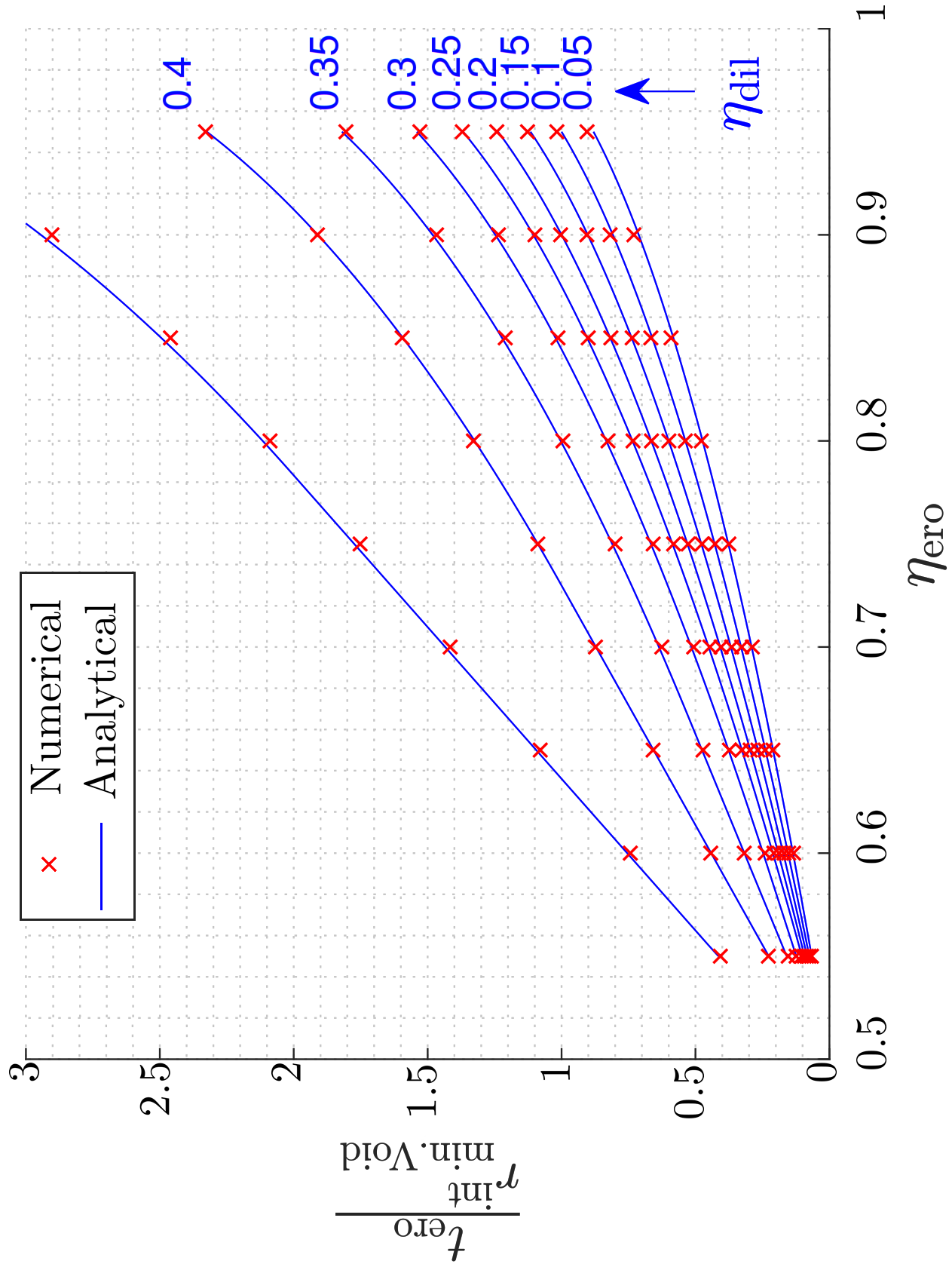}
		\caption{}
		\label{fig:fig_1_d}
	\end{subfigure}
	\vspace{-222mm}\\
	\hspace{-80mm}\Large{(a)} \hspace{77mm} \Large{(b)}
	\vspace{106mm}\\
    \hspace{-80mm}\Large{(c)} \hspace{77mm} \Large{(d)}
	\vspace{98mm}
	\caption{Graphical relationships between the minimum length scale and $\eta_\mathrm{ero}$, $\eta_\mathrm{int}$, $\eta_\mathrm{dil}$, and $r_\mathrm{fil}$. The minimum length scale is defined for the intermediate design. The graphs are built using both, the analytical and the numerical method. Graphs (a) and (b) are designed to obtain the projection thresholds and the filter radius, respectively. Graphs (c) and (d) are designed to obtain the dilation and erosion distances, respectively.}
	\label{fig:fig_1}
\end{figure*}

Firstly, a number of advantages have been observed when defining the intermediate design with a threshold ${\Threshold{int}} = 0.5$. For instance, the projection features lower amounts of intermediate densities compared to those projections that use a threshold other than 0.5 \citep{Xu2010,Wang2011,DASILVA2019_2}, and a threshold ${\Threshold{int}}$ set to 0.5 provides the same size ranges (0.5) for the erosion and dilatation thresholds, which is convenient for reducing rounding errors, since small differences in projection thresholds could be insensitive to the minimum size when using a coarse discretization of the design domain \citep{Qian2013}. Secondly, for a particular combination of thresholds, the filter size ($r_\mathrm{fil}$) can become considerably larger than the desired minimum length scale, which could significantly increase computational requirements \citep{Lazarov2011}. This can be seen in Figs.~\ref{fig:MinVoidFilter} and \ref{fig:fig_1_b}. These figures show graphs that relate the filter size ($r_\mathrm{fil}$) to the minimum size of solid or void phase and to the erosion or dilation threshold. These graphs show that the closer ${\Threshold{ero}}$ and ${\Threshold{dil}}$ are to ${\Threshold{int}}$, the larger the filter radius, which inevitably increases computational requirements. Under this observation, we recommend choosing ${\Threshold{ero}} \geq 0.75$ and ${\Threshold{dil}}\leq 0.25$, thus it is ensured that $r_\mathrm{fil}\leq 2 r_\mathrm{min.Solid}^\mathrm{int}$ and $r_\mathrm{fil}\leq 2 r_\mathrm{min.Void}^\mathrm{int}$. Thirdly, erosion and dilation thresholds too distant or too close to the intermediate threshold increases oscillations of design variables during the optimization process. In general, a good compromise is to choose $0.10 \leq {\Threshold{dil}} \leq 0.4$ and $0.60 \leq {\Threshold{ero}} \leq 0.9$.

The first observation removes an unknown from the system of equations, since ${\Threshold{int}}$ is set to 0.5, while the second and third observations limit the range of the erosion and dilation thresholds. Thus, for a user-defined minimum length scale, it is possible to develop an algorithm that solves the system of equations considering the three observations. Here we propose an algorithm based on graphic relationships, so that the reader can easily find the desired parameters without the need to resort to a computational algorithm. Nonetheless, we also provide as supplementary material a code written in MATLAB named \texttt{SizeSolution.m} that performs the procedure described below. It is important to note that the above observations are based on numerical tests considered for specific optimization problems formulated in the density approach. Therefore, it is possible that under other topology optimization approaches or formulations the above observations are no longer valid. However, the proposed procedure can be applied for any other value of ${\Threshold{int}}$, or any other combination of parameters that the user may consider convenient.

\begin{figure}[t!]
	\captionsetup{width=1.00\linewidth}
    	\centering
		\includegraphics[width=0.95\linewidth]{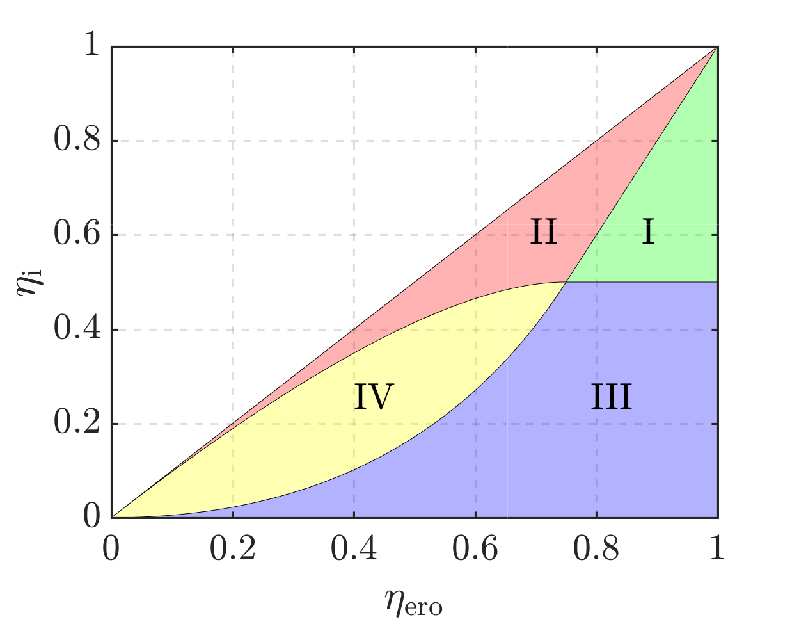}
	\caption{Graphic representation of the applicability of the set of equations in Table \ref{tab:SolidZoneCondition_v2} proposed for the solid phase. The denomination of the zone corresponds to the row number of Table \ref{tab:SolidZoneCondition_v2}.}
	\label{fig:SolidZones}
\end{figure}

Given that the ranges of application of the equations in Table \ref{tab:SolidZoneCondition_v2} are defined as a function of the projection thresholds ($\eta_i$, ${\Threshold{ero}}$ and ${\Threshold{dil}}$), it is rather simple to construct graphs with respect to them. For instance, as shown in Fig.~\ref{fig:SolidZones}, to find the range of application of the equations defining the minimum size of the solid phase, it is simply necessary to know in which region of Fig.~\ref{fig:SolidZones} the projection $\eta_i$ falls.

The first graph proposed in this work is shown in Fig.~\ref{fig:fig_1_a} and gathers 4 parameters, the minimum length scale ($r_\mathrm{min.Solid}^\mathrm{int}$ and $r_\mathrm{min.Void}^\mathrm{int}$) and the projection thresholds (${\Threshold{ero}}$ and ${\Threshold{dil}}$), provided that ${\Threshold{int}} = 0.5$. In this graph, the user can easily find the set of erosion and dilation threshold that leads to the desired length scale. Then, the user can access the graph in Fig.~\ref{fig:fig_1_b} to obtain the filter radius. For example, for the following minimum length scale, $r_\mathrm{min.Solid}^\mathrm{int}=3$ elements and $r_\mathrm{min.Solid}^\mathrm{int}=3$ elements, the graph in Fig.~\ref{fig:fig_1_a} is accessed with a value of 1.0 for the ordinate. According to the aforementioned observations, the combination of thresholds [${\Threshold{ero}}$, ${\Threshold{dil}}$] that can be chosen among others are:
\begin{subequations} \label{Eq:Threshold_rMinSolid_rMinVoid}
	\begin{align}
	[{\Threshold{ero}}, {\Threshold{dil}}] &= [0.75,\; 0.25]  \; ,\\[1.0ex] 
	[{\Threshold{ero}}, {\Threshold{dil}}] &= [0.80,\; 0.20]  \; ,\\[1.0ex] 
	[{\Threshold{ero}}, {\Threshold{dil}}] &= [0.85,\; 0.15]  \; ,\\[1.0ex] 
	[{\Threshold{ero}}, {\Threshold{dil}}] &= [0.90,\; 0.10]  \; .
	\end{align}
\end{subequations}

Arbitrarily, [0.75, 0.25] is selected, and from the graph in Fig.~\ref{fig:fig_1_b}, it is obtained that $r_\mathrm{fil}=2r_\mathrm{min.Solid}^\mathrm{int} = 6$ elements.

It should be noted that the graph in Fig.~\ref{fig:fig_1_a} is constructed considering a resolution of 0.05 in the projection thresholds. This is due to the fact that the discretization of the filter radius in topology optimization is generally coarse, and decimal numbers smaller than 0.05 in the threshold value have usually a negligible effect on the minimum length scale of the optimized design. However, if a better resolution is required, the user can resort to the attached code \texttt{SizeSolution.m}.

% ================ FIGURE THAT GATHERS ALL THE GRAPHS =============
% =================================================================
% ================================================================
% ================================================================

As previously mentioned, the erosion and dilation distances could be required, for instance, when implementing maximum size restrictions \citep{Fernandez2020}. These distances can be easily obtained from the equations of Table \ref{tab:SolidZoneCondition_v2}. As shown in Fig.~\ref{fig:fig_1_v2_d}, the dilation distance, denoted by $t_\mathrm{dil}$, is the offset distance between the intermediate and dilated designs. Analogously, the erosion distance, denoted by $t_\mathrm{ero}$, is the offset distance between the intermediate and eroded designs. Therefore:

\begin{equation} \label{Eq:Dilation_and_Erosion}
\begin{matrix}
t_\mathrm{dil} &= r_\mathrm{min.Solid}^\mathrm{dil} {(\Threshold{dil},\Threshold{ero},r_\mathrm{fil})} - r_\mathrm{min.Solid}^\mathrm{int} {(\Threshold{int},\Threshold{ero},r_\mathrm{fil})}
\vspace{2mm}\\
t_\mathrm{ero} &= r_\mathrm{min.Void}^\mathrm{ero} {(\Threshold{ero},\Threshold{dil},r_\mathrm{fil})} - r_\mathrm{min.Void}^\mathrm{int} {(\Threshold{int},\Threshold{dil},r_\mathrm{fil})}
\end{matrix}
\end{equation}
\noindent where the minimum size of the solid phase in the dilated design ($r_\mathrm{min.Solid}^\mathrm{dil}$) is obtained by choosing $\eta_i = {\Threshold{dil}}$ in the equations of Table \ref{tab:SolidZoneCondition_v2}. Analogously, the minimum size of the void phase in the eroded design ($r_\mathrm{min.Void}^\mathrm{ero}$) is obtained by choosing $\eta_i = {\Threshold{ero}}$. In this way, analytical expressions can be obtained for the dilation and erosion distances. 

It is noted from Eq.~\eqref{Eq:Dilation_and_Erosion} that the erosion and dilation depend on the 3 projection thresholds and on the size of the filter, namely, $t_\mathrm{ero}(r_\mathrm{fil},{\Threshold{ero}},{\Threshold{int}},{\Threshold{dil}})$ and $t_\mathrm{dil}(r_\mathrm{fil},{\Threshold{ero}},{\Threshold{int}},{\Threshold{dil}})$. Therefore, if for some particular reason the user needs to impose the minimum length scale, and the erosion and dilation distances, then the following system of equations must be solved: 
\begin{equation} \label{Eq:Determined_System}
\begin{split}
r_\mathrm{min.Solid}^\mathrm{int} & = r_\mathrm{min.Solid}^\mathrm{int}(r_\mathrm{fil},{\Threshold{int}},{\Threshold{ero}}) \vspace{1mm}\\
r_\mathrm{min.Void}^\mathrm{int}  & = r_\mathrm{min.Solid}^\mathrm{int}(r_\mathrm{fil},{\Threshold{int}},{\Threshold{dil}}) \vspace{2mm}\\
t_\mathrm{ero} & =  r_\mathrm{min.Solid}^\mathrm{int}(r_\mathrm{fil},{\Threshold{ero}},{\Threshold{int}},{\Threshold{dil}}) \vspace{1mm}\\
t_\mathrm{dil} & =  r_\mathrm{min.Solid}^\mathrm{int}(r_\mathrm{fil},{\Threshold{ero}},{\Threshold{int}},{\Threshold{dil}}) \\
\end{split}
\end{equation} 

The system in Eq.~\eqref{Eq:Determined_System} is determined and there is only one combination of parameters [$r_\mathrm{fil}$, ${\Threshold{ero}}$, ${\Threshold{int}}$, ${\Threshold{dil}}$] that leads to the desired length scale [$r_\mathrm{min.Solid}^\mathrm{int}$, $r_\mathrm{min.Void}^\mathrm{int}$, $t_\mathrm{ero}$, $t_\mathrm{dil}$].  This can be illustrated with the following example. For each combination of thresholds [$\Threshold{ero}$, $\Threshold{dil}$] given in Eq.~\eqref{Eq:Threshold_rMinSolid_rMinVoid}, the erosion and dilation distances [${t_\mathrm{ero}}$ , ${t_\mathrm{dil}}$ ] are provided, as follows:
\begin{subequations} \label{Eq:Threshold_rMinSolid_rMinVoid_offset}
	\begin{align}
	[{\Threshold{ero}}, {\Threshold{dil}}] &= [0.75,\; 0.25] \; ,
	[{t_\mathrm{ero}} , {t_\mathrm{dil}} ]  = [1.76,\; 1.76] \\[0.5ex] 
	[{\Threshold{ero}}, {\Threshold{dil}}] &= [0.80,\; 0.20]  \; ,
	[{t_\mathrm{ero}} , {t_\mathrm{dil}} ]  = [1.99,\; 1.99] \\[0.5ex] 
	[{\Threshold{ero}}, {\Threshold{dil}}] &= [0.85,\; 0.15]  \; ,
	[{t_\mathrm{ero}} , {t_\mathrm{dil}} ]  = [2.21,\; 2.21] \\[0.5ex] 
	[{\Threshold{ero}}, {\Threshold{dil}}] &= [0.90,\; 0.10]  \; ,
	[{t_\mathrm{ero}} , {t_\mathrm{dil}} ]  = [2.43,\; 2.43]  
	\end{align}
\end{subequations}

It is recalled that all combinations of thresholds [$\Threshold{ero}$, $\Threshold{dil}$] in Eq.~\eqref{Eq:Threshold_rMinSolid_rMinVoid_offset} lead to the minimum length scale $r_\mathrm{min.Solid}^\mathrm{int} = r_\mathrm{min.Void}^\mathrm{int} = 3$ elements. However, they all result in different erosion and dilation distances.

To the authors' knowledge, the need to impose specific values for the erosion and dilation distances has not yet been claimed, so this manuscript is limited to providing these values for a given combination of parameters. The proposed graphs are shown in Fig.~\ref{fig:fig_1_c} and \ref{fig:fig_1_d}, which depend on the erosion and dilation thresholds and on the corresponding minimum size. For instance, considering $[{\Threshold{ero}}, {\Threshold{dil}}]$ = [0.75, 0.25], the erosion and dilation distances are $t_\mathrm{ero} = 0.58 \times 3 $ elements and $t_\mathrm{dil} = 0.58 \times 3 $ elements, where 0.58 is obtained from the graphs.

\section{Sources of error} \label{sec:5}

The analytical procedure developed in this article allows to quickly obtain the set of filtering and projection parameters that imposes the desired minimum length scale. However, in practice, the length scale of the optimized design often differs from the desired values \citep{Qian2013}. This is mainly due to the fact that the analytical method assumes (i) a continuous domain, (ii) a perfect Heaviside projection ($\beta \to \infty$), and (iii) that the eroded and dilated fields project an infinitesimal minimum size in the solid and void phases, respectively. The assumptions (i), (ii), and (iii) are not met in topology optimization mainly because the design domains are discretized into finite elements and because the Heaviside projection is smoothed. To assess the error introduced by these assumptions, in the following we compare the analytical method with the numerical method proposed by \citet{Wang2011}, which considers a discrete domain and a smoothed Heaviside function.

The procedure for obtaining the minimum size using the numerical method is analogous to the analytical method but now using a one-dimensional domain discretized into $N$ elements. For example, to obtain the minimum size of the solid phase the three-field scheme is applied to a one-dimensional design domain $\bm{\rho}$ containing solid elements in a length $h$. Then, the size of the eroded, intermediate and dilated fields are measured, and the length $h$ is adjusted so that the resulting eroded field has an infinitesimal size (1 solid element). This process is repeated several times for different values of ${\Threshold{ero}}$, and the resulting minimum size is normalized with respect to the size of the chosen filter. For implementation details regarding the numerical method, the reader is referred to the works of \citet{Wang2011} and \cite{Fernandez2020}, and to the attached MATLAB code named \texttt{NumericalSolution.m}.

To validate the analytical and numerical methods, the latter is implemented using a Heaviside function at $\beta=500$, a filter radius of 1000 elements and a design domain discretized into 10 thousand elements, so the numerical method can be considered continuous and under similar assumptions than the analytical one. The relationships obtained with the numerical method can be seen in the graphs of Fig.~\ref{fig:fig_1}. The agreement of results from both methods allows us to validate the set of equations provided in Table \ref{tab:SolidZoneCondition_v2}. 

The three sources of error that affect the analytical method are discussed below.

\subsection{Continuous design domain}\label{Sec:rounding_error}

\begin{figure}[t!]
\begin{subfigure}{0.92\linewidth}
    	\centering
		\includegraphics[width=0.9\linewidth]{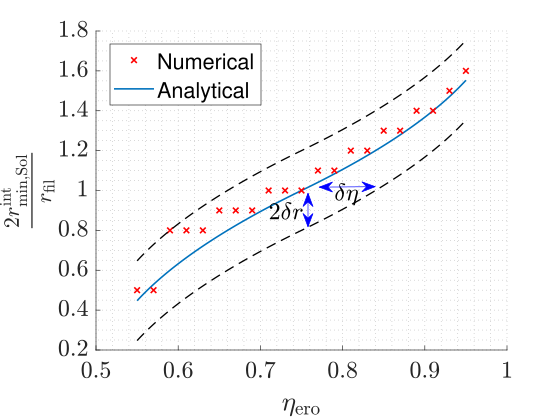}
		\caption{$N = 100$, $r_\mathrm{fil}=10$.}
		\label{fig:RoundingError10}
	\end{subfigure}
	\vspace{3mm}\\
	\begin{subfigure}{0.92\linewidth}
    	\centering
		\includegraphics[width=0.9\linewidth]{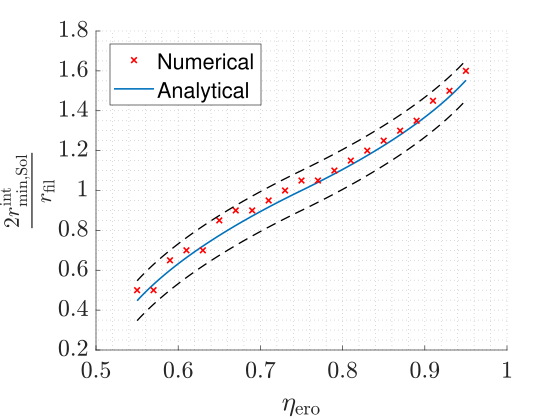}
		\caption{$N = 200$, $r_\mathrm{fil}=20$}
		\label{fig:RoundingError20}
	\end{subfigure}
\caption{Relationship between the minimum size of the solid phase, the eroded threshold, and the filter radius. The dotted line represents the analytical relationship including the rounding error.}
\label{fig:RoundingError}
\end{figure}

When the design domain is discretized using finite elements, the radii that define the minimum sizes ($r_\mathrm{fil}$, $\rminsolid{int}$, $\rminsolid{dil}$, $\rminvoid{int}$ and $\rminvoid{ero}$) are defined by a discrete number of elements. In this case, the rounding error is $\pm$ 1 finite element in the radius that defines the minimum size. Therefore, to reduce this error, it is sufficient to reduce the size of the elements by mesh refinement. To illustrate this remark graphically, we consider the following two discretizations of a one-dimensional domain, one containing 100 elements and the other containing 200 elements. In our implementation, the size of the filter is chosen equal to $10\%$ of the domain size. Therefore, for the discretizations containing 100 and 200 elements, the filter radius contains 10 and 20 elements, respectively. For each discretization, the relationship between the minimum size of the solid phase, the size of the filter and the erosion threshold is plotted in Fig.~\ref{fig:RoundingError}, using both the analytical and the numerical method. 

As the analytical method does not depend on the discretization, it provides identical relationships on both discretizations. However, the rounding error associated to the analytical method does depend on the discretization and can be plotted as an offset of the analytical curve, as shown by the dashed curves in Fig.~\ref{fig:RoundingError}. The rounding error, denoted as $\delta r$, corresponds to 1 element in the radius, i.e.~2 elements in the diameter. That is, the vertical offset of the analytical curve due to the rounding error is equal to $2\delta r=2/10$ and to $2\delta r=2/20$ in Figs.~\ref{fig:RoundingError10} and \ref{fig:RoundingError20}, respectively.

The error curves agree with the results obtained from the numerical method, which shows that the error of the analytical method coming from the assumption of a continuous domain can be easily estimated. This estimation allows to know either the margin of error of the minimum size ($\delta r$) for a given analytical threshold (${\Threshold{ero}}$ or ${\Threshold{dil}}$), or the margin of error of the analytical threshold ($\delta \eta$) for a given minimum size, as shown in Fig.~\ref{fig:RoundingError10}.

If the user of the robust approach needs to impose precisely the minimum size, he can choose the projection thresholds (${\Threshold{ero}}$ and ${\Threshold{dil}}$) that lead to a larger filter radius, or he can use the numerical method \citep{Wang2011} using a discretization representative of the design domain to be optimized.

\subsection{Perfect Heaviside projection}

The analytical method developed in this work assumes a perfect Heaviside function, which results in a projected field $\bm{\bar{\rho}}$ containing discrete densities (0 and 1). However, in topology optimization, the Heaviside function is smoothed resulting in regions of intermediate densities that lie on the surface of the optimized structure. For this reason, in practice, once the optimized solution ${\ProjField{int}}$ is obtained, a cut-off value $\varepsilon$ is usually added to the projected densities in order to get the optimized structure. Therefore, the final design intended for manufacturing is obtained by projecting the design ${\ProjField{int}}$ as follows:

\begin{equation}
 {\bar{\rho}^{\varepsilon}}_i = 
 \left\lbrace
   \begin{matrix}
     1 \;\;, & \text{if} \;\;\;\; \bar{\rho}_{\mathrm{int}(i)} \geq \varepsilon \vspace{1mm}\\
     0 \;\;, & \text{otherwise}
   \end{matrix}
 \right. 
\end{equation} 

Since the post-processed design $\bm{\bar{\rho}^{\varepsilon}}$ is the one intended for manufacturing, the analytical expressions relating the minimum length scale and the parameters that define it must be elaborated for the field $\bm{\bar{\rho}^{\varepsilon}}$. As will be explained hereafter, the equations provided in Table \ref{tab:SolidZoneCondition_v2} can be easily adapted to consider the cut-off value $\varepsilon$.

\begin{figure}[t!]
	\captionsetup{width=1.00\linewidth}
    	\centering
		\includegraphics[width=0.85\linewidth]{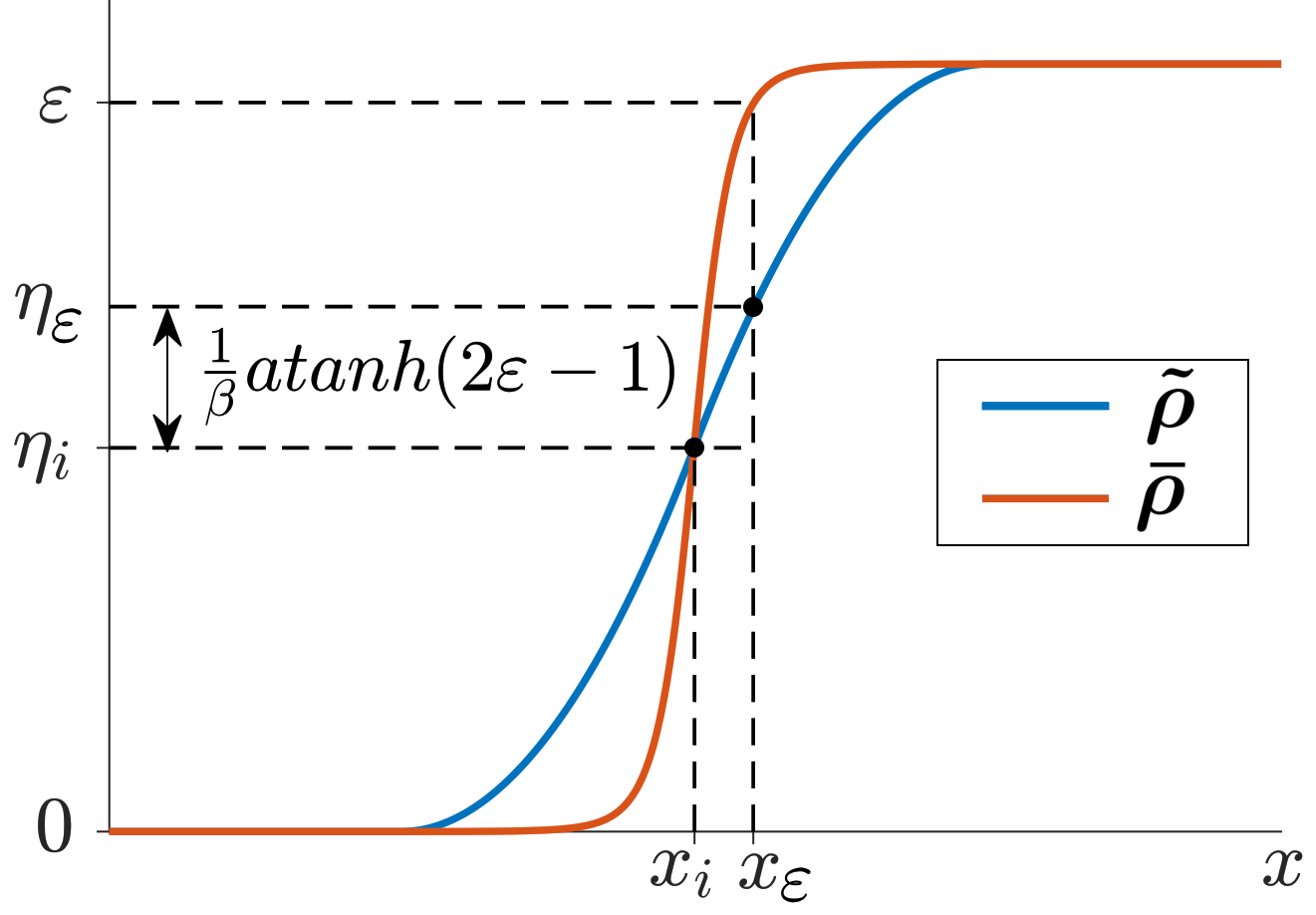}
	\caption{{ Illustration of a transition region where projected densities reach intermediate values.}}
	\label{fig:betaError}
\end{figure}

Consider Fig.~\ref{fig:betaError} illustrating a transition zone between the solid structure and the void region. In blue it is shown the filtered field $\bm{\tilde{\rho}}$ and in red the projected field $\bm{\bar{\rho}}$ obtained with $\beta=32$. In this illustration, the density cut-off is defined as $\varepsilon=0.95$. A perfect Heaviside projection ($\beta \to \infty$) would define the size of the solid zone at the coordinate $x_i$. However, the smoothed Heaviside projection defines the solid/void transition in the coordinate $x_\varepsilon$. The idea then is to find the projection threshold $\eta_\varepsilon$ for which a perfect Heaviside projection produces the same solid/void transition coordinate than the one that is obtained with a smoothed Heaviside projection with a threshold $\eta_i$. Thus, the value of the filtered density present at the coordinate $x_\varepsilon$ must be found. To this end, from Fig.~\ref{fig:betaError} it is observed that:
\begin{equation}
\varepsilon = \frac{\tanh{(\beta \eta_i)} \tanh{(\beta(\eta_\varepsilon -\eta_i))}}
{\tanh{(\beta \eta_i)} \tanh{(\beta(1 - \eta_i))}}
\label{Eq:Heaviside_varepsilon}
\end{equation}

Assuming $\beta>10$ at the end of the optimization process, which is a common practice when dealing with Heaviside projection, Eq.~\eqref{Eq:Heaviside_varepsilon} can be simplified to:
\begin{equation}
\tanh{(\beta(\eta_\varepsilon -\eta_i))} = 2\varepsilon -1 
\label{Eq:Heaviside_varepsilon_2}
\end{equation}

From Eq.~\eqref{Eq:Heaviside_varepsilon_2}, the filtered density $\eta_\varepsilon$ can be obtained: 
\begin{equation}
\eta_\varepsilon = \eta_i + \frac{1}{\beta} \mathrm{atanh}(2\varepsilon -1) 
\label{Eq:Heaviside_varepsilon_3}
\end{equation}

We recall that the filtered density $\eta_\varepsilon$ represents the value for which a perfect Heaviside function with threshold $\eta_\varepsilon$ and a smoothed Heaviside function with threshold $\eta_i$ result in the same size for the solid and void phases. Therefore, to take into account the intermediate densities resulting from a smoothed Heaviside function, the shifting term $\mathrm{atanh}(2\varepsilon -1)/\beta$ must be added to the thresholds involved in the analytical equations (Table \ref{tab:SolidZoneCondition_v2}). Note that if $\beta \to \infty$ or $\varepsilon=0.50$, then the shifting term is zero and $\eta_\varepsilon = \eta_i$.

\begin{figure}[t!]
	\captionsetup{width=1.00\linewidth}
    	\centering
    	\begin{subfigure}{1.\linewidth}
		\includegraphics[width=0.95\linewidth]{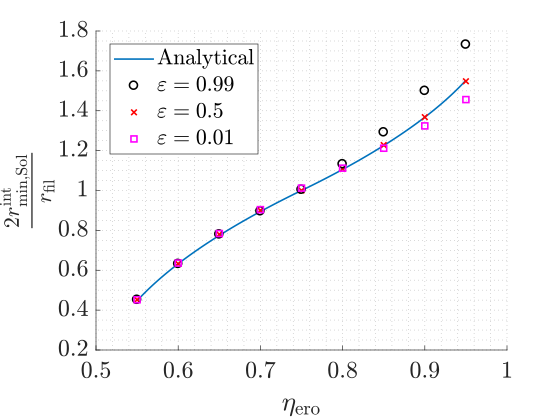}
		\caption{}
		\label{fig:betaErrorEx2}
		\end{subfigure}
		\\
		\begin{subfigure}{1.\linewidth}
		\includegraphics[width=0.95\linewidth]{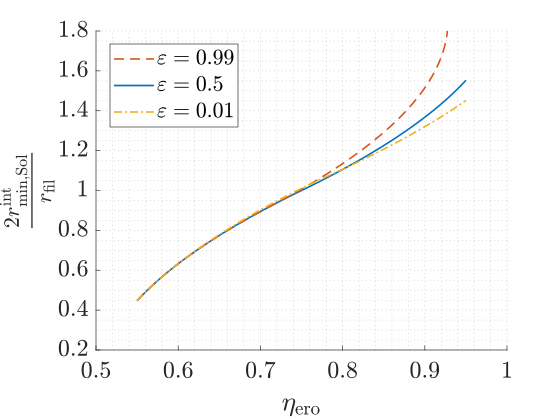}
		\caption{}
		\label{fig:betaErrorEx}
		\end{subfigure}
	\caption{(a) Minimum length scale obtained with the analytical method assuming a perfect Heaviside projection ($\beta \to \infty$), and with the numerical method using a smoothed Heaviside projection ($\beta=32$) and three different cut-off values $\varepsilon$. (b) The analytical curves obtained by adding the shifting term $\mathrm{atanh}(2\varepsilon -1)/\beta$ to the thresholds.}
	\label{fig:betaErrorEx22}
\end{figure}

To assess the error introduced by a smoothed Heaviside projection, the minimum length scale obtained with the analytical and the numerical methods are compared. The analytical method is applied with $\beta \to \infty$, so it does not consider the shifting term on the thresholds. For the numerical method we use $\beta=30$, a filter radius equal to 200 and a one-dimensional domain discretized into 2000 elements in order to avoid rounding errors and isolate the effect of using a smoothed Heaviside projection. In the numerical method, three cut-off values are used, $\varepsilon=0.01$, $\varepsilon=0.50$ and $\varepsilon=0.99$. The results are shown in Fig.~\ref{fig:betaErrorEx2}.

The results show that the influence of using a smoothed Heaviside function is observed when $\varepsilon \neq 0.5$ and when ${\Threshold{ero}}>0.75$. For this reason, we recommend to use a density cut-off value equal to 0.5, as this avoids the need to add the shifting term to the projection thresholds. Nevertheless, the user can easily adjust the analytical curves for $\varepsilon \neq 0.5$, since it is only required to add $\mathrm{atanh}(2\varepsilon -1)/\beta$ to the projection thresholds, as shown in the corrected curves in Fig.~\ref{fig:betaErrorEx}.

It is worth mentioning that probably the effect of using a smoothed Heaviside function is only seen for ${\Threshold{ero}} \geq 0.75$ because the projected field becomes more discrete as the projection threshold is closer to 0.5 \citep{Wang2011}.

\subsection{Infinitesimal size} \label{Sec:Infinitesimal_Size}

It is recalled that the analytical and numerical methods are developed under the condition of robustness which ensures that the manufactured design will feature a good performance even if the blueprint design is uniformly thinned or thickened during the manufacturing process. This condition dictates that the structural members/cavities must be projected by at least one solid/void element in the case of an eroded/dilated design.

\begin{figure}[t!]
    	\centering
		\includegraphics[width=0.85\linewidth]{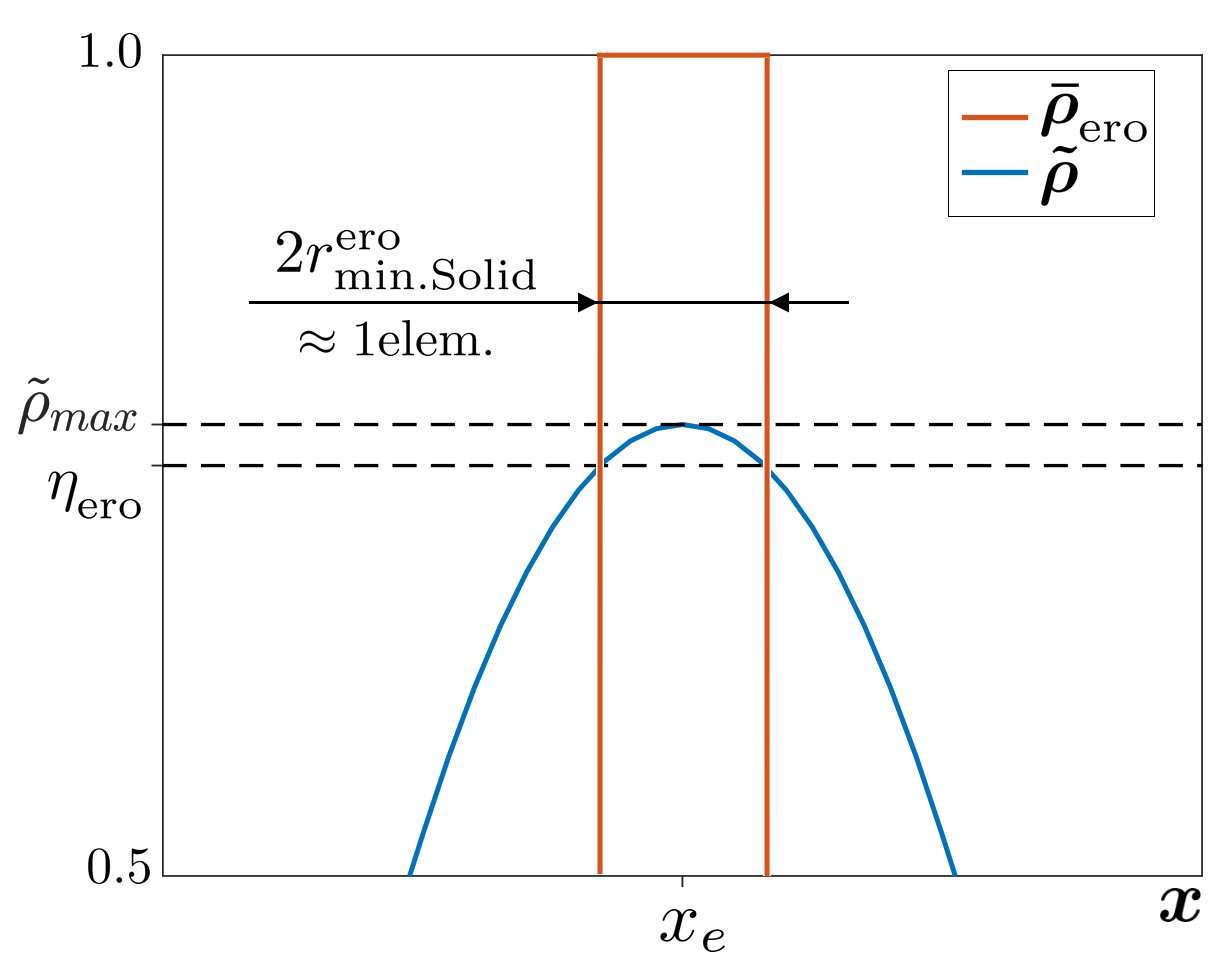}
		\caption{Illustration in which an infinitesimal size is not reached for the eroded projection.}
\label{fig:AvoidError}
\end{figure}

To obtain the minimum size, the analytical method assumes that the size of the structural members/cavities is infinitesimal ($\approx 0$) in the eroded/dilated design, but in practice this does not occur. In topology optimization, it is observed that the smallest structural members/cavities in the eroded/dilated design are composed of one or a few elements, generally described with intermediate densities. Thus, in practice, the minimum size of the solid phase in the eroded design ($r_\mathrm{min.Solid}^\mathrm{ero}$) results in a discrete number of elements. This can be seen in Fig.~\ref{fig:AvoidError}, which is build using the numerical method. There, a filtered field and its eroded projection under the robust condition are shown. It can be seen that the minimum size of the eroded design is not infinitesimal, and therefore the minimum size of the solid phase in the intermediate and dilated design would be larger than the value predicted by the analytical method.

\begin{figure}[t!]
    	\centering
		\includegraphics[width=0.95\linewidth]{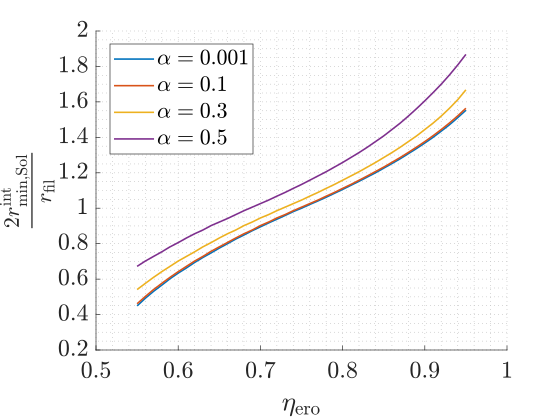}
		\caption{Effect of not reaching an infinitesimal size in the condition of robustness $r_\mathrm{min.Solid}^\mathrm{ero}=0$. This condition is assessed by imposing $r_\mathrm{min.Solid}^\mathrm{ero}= \alpha \: r_\mathrm{fil}$.}
\label{fig:Infinitesimal_error}
\end{figure}

The minimum size of the solid and void phases in the eroded and dilated designs ($r_\mathrm{min.Solid}^\mathrm{ero}$ and $r_\mathrm{min.Void}^\mathrm{dil}$) are defined by the size of the elements that discretize the design space. In addition, the amount of elements with intermediate densities defining the minimum size depends on the steepness of the smoothed Heaviside function ($\beta$), therefore, the error introduced by assuming an infinitesimal size in the robust condition is related to the two sources of error mentioned previously. To isolate the effect of the infinitesimal size in the robust condition and illustrate the error introduced by the assumption of infinitesimal size, we make use of the numerical method. The numerical method is implemented in a discretized domain containing a large number of elements ($10^4$ elements) and using a large filter size ($10^3$ elements) to simulate a continuous domain. The steepness parameter of the smoothed Heaviside projection is set as $\beta=512$. Thus, the size and density of the elements are excluded from the analysis. The effect of not achieving an infinitesimal size in the condition of robustness is intentionally imposed in the numerical code. For this, the robustness condition is considered satisfied if $r_\mathrm{min.Solid}^\mathrm{ero} = \alpha \: r_\mathrm{fil}$. Considering that $r_\mathrm{min.Solid}^\mathrm{ero}$ represents one finite element in topology optimization, and that representative values for $r_\mathrm{fil}$ are between 2 and 10 elements, it is reasonable to consider values for $\alpha$ ranging from 0.1 to 0.5. Taking into consideration the above, we build the graphical solutions that relate the minimum size in the solid phase ($r_\mathrm{min.Solid}^\mathrm{int}$), the filter radius and the erosion threshold. The graph is shown in Fig.~\ref{fig:Infinitesimal_error}.

The graph shows that the error of not achieving an infinitesimal size in the eroded design produces a minimum size of the solid phase bigger than the value predicted by the analytical method. The error related to the infinitesimal size is low in comparison to the rounding error, so the latter would be the most relevant source of error from the analytical method that assumes a continuous design domain.

In summary, this section discussed the scope of the analytical method \citep{Qian2013} by comparing it with the numerical one \citep{Wang2011}, both developed for a one-dimensional design domain. To this end, different sources of error were examined. In general, the errors can be controlled by mesh refinement or by correcting the projection thresholds according to the cut-off $\varepsilon$ value. In the following section, the analytical method is assessed using 2D-topology optimization problems. 

\section{Numerical examples and discussion}\label{sec:6}

This section examines the reliability of the analytical expressions provided in Table \ref{tab:SolidZoneCondition_v2}. To this end, a set of 2D topology optimization problems are solved, from which the length scale is measured graphically and compared with the imposed values. Then, some designs obtained with maximum size constraints are provided to illustrate the use of the erosion and dilation distances. Finally, this section provides a remark regarding the simplified robust formulation, where the intermediate and dilated designs are removed from the objective function.

\subsection{Minimum length scale}

\begin{figure}[t!]
    	\centering
		\includegraphics[width=0.65\linewidth]{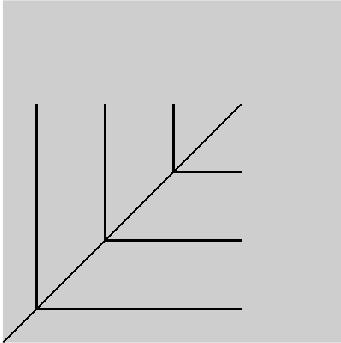}
		\caption{Initial distribution of design variables considered for the thermal compliance minimization problem.}
\label{fig:Skeleton}
\end{figure}

The minimum length scale is assessed using the heat exchanger design problem described in Section \ref{sec:2}. A set of results is generated from this design problem, which differ in the desired minimum length scale and in the discretization used. Specifically, three sets of solutions are obtained by discretizing the design domain into $100 \times 100$, $200 \times 200$ and $400 \times 400$ quadrilateral elements. In addition, three different length scales are prescribed for each discretization, which are reported as the ratio between the minimum size of the solid phase and the minimum size of the void phase, i.e.~$r_\mathrm{min.Solid}^\mathrm{int}/r_\mathrm{min.Volid}^\mathrm{int}$. The chosen ratios are $1/2$, $1/1$ and $2/1$. The minimum size in the solid phase is the same in all scenarios and is defined as a physical dimension. In number of finite elements, the radius that defines the minimum size of the solid phase ($r_\mathrm{min}^\mathrm{int}$) is equal to 1, 2 and 4, for the discretizations that use $100^2$, $200^2$ and $400^2$ elements, respectively. It is well known that the initial values of design variables have a huge influence on the resulting topology when it comes to thermal compliance minimization \citep{Yan2018}, hence, in order to facilitate the comparison of results, we impose a base structure as a starting point, which is shown in Fig.~\ref{fig:Skeleton}. Before presenting the results, the procedure to obtain the minimum length scale from the optimized designs is detailed. 

%===========================
%===========================

\begin{figure}[t!]
    \captionsetup[subfigure]{labelformat=empty}
    	\centering
    		\begin{subfigure}{0.95\linewidth}
		\centering \includegraphics[width=0.8\linewidth]{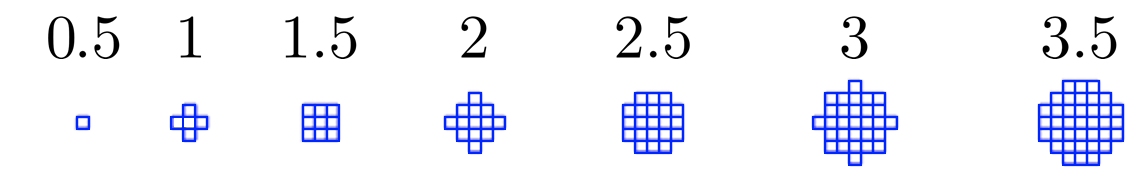}
		\end{subfigure}
		\\
    		\begin{subfigure}{0.95\linewidth}
    		\centering \includegraphics[width=0.95\linewidth]{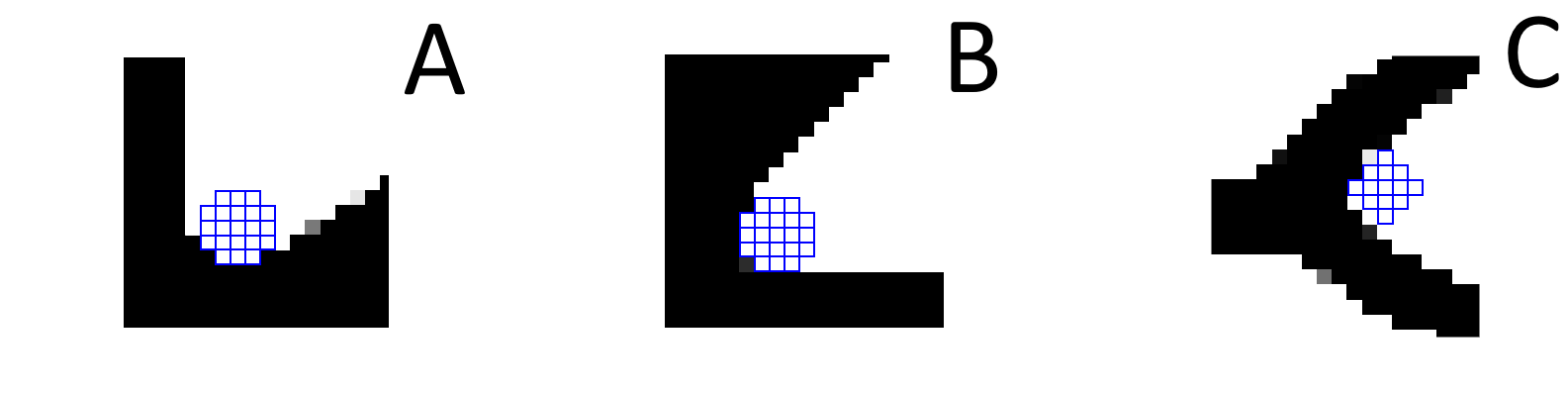}
		\end{subfigure}
		\\
    		\begin{subfigure}{0.95\linewidth}
    		\centering \includegraphics[width=0.75\linewidth]{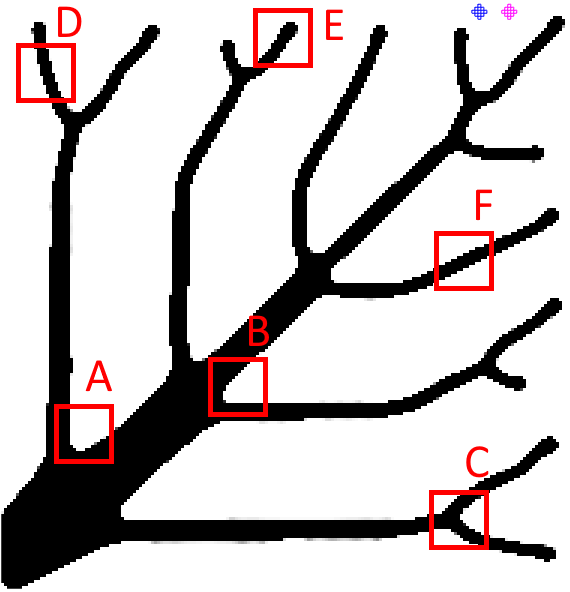}
		\end{subfigure}
		\\
		\begin{subfigure}{0.95\linewidth}
    		\centering \includegraphics[width=0.85\linewidth]{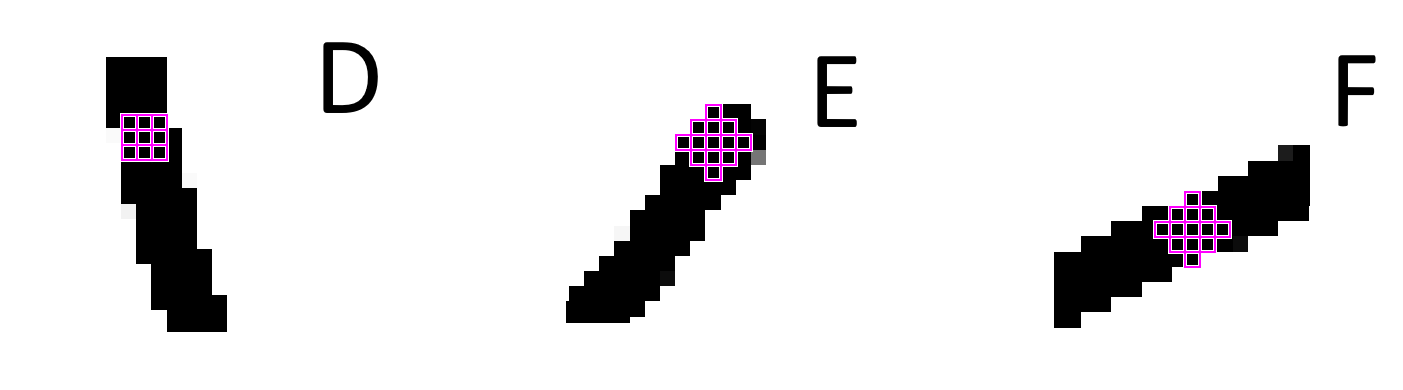}
		\end{subfigure}
		\\
		\begin{minipage}{0.01\textwidth}
			\subcaption{\label{fig:IllustrationOfMeasureSize_a}}
		\end{minipage}	
		~
		\begin{minipage}{0.01\textwidth}
			\subcaption{\label{fig:IllustrationOfMeasureSize_b}}
		\end{minipage}
		~
		\begin{minipage}{0.01\textwidth}
			\subcaption{\label{fig:IllustrationOfMeasureSize_c}}
		\end{minipage}
		~
		\begin{minipage}{0.01\textwidth}
			\subcaption{\label{fig:IllustrationOfMeasureSize_d}}
		\end{minipage}
		\vspace{-117mm}\\
		\hspace{-78mm}\footnotesize{(a)}
		\vspace{13mm}\\
		\hspace{-79mm}\footnotesize{(b)}
		\vspace{37mm}\\
		\hspace{-79mm}\footnotesize{(c)}
		\vspace{38mm}\\
		\hspace{-79mm}\footnotesize{(d)}
		\vspace{10mm}\\
		\vspace{-2mm}
		\caption{Illustration of the minimum length scale measurement. In (a) the test regions and their radius sizes given in number of finite elements. In (b) the minimum size of the void phase, (c) the optimized heat exchanger, and (d) the minimum size of the solid phase. The imposed length scales are graphically shown in (c), at the upper left corner of the design.}
		\label{Fig:IllustrationOfMeasureSize}
\end{figure}

The minimum size of the solid phase is measured graphically by counting the number of finite elements that define the size of the thinnest structural branch. Similarly, the minimum size of the void phase is measured by counting the elements in the radius of the largest circumference that can be inscribed at the re-entrant corners of the design. For example, consider the design of Fig.~\ref{fig:IllustrationOfMeasureSize_c} where the minimum length scales $r_\mathrm{min.Solid}^\mathrm{int} = \mathrm{min.Void}^\mathrm{int} =  2$ elements are imposed. To determine the real minimum size of the void phase, the largest circle in Fig.~\ref{fig:IllustrationOfMeasureSize_a} that falls into the re-entrant corners of the design is identified. The corners analyzed are those that form a sharp angle between two structural branches. Fig.~\ref{fig:IllustrationOfMeasureSize_b} shows three representative re-entrant corners of the design depicted in Fig.~\ref{fig:IllustrationOfMeasureSize_c}. From there it is observed that the minimum size is given by a circle of radius 2 elements (zone C). Similarly, the largest region that fits into the thinnest structural members is determined, as shown in Fig.~\ref{fig:IllustrationOfMeasureSize_d}. There, the minimum size of the solid phase is given by a circle of radius 1.5 elements (zone D).

After describing the test case, the obtained results are presented. Table \ref{Tab:9_HeatExchangers} contains the nine results generated in this example (3 length scales $\times$ 3 discretizations). The imposed minimum length scales are reported graphically next to each solution. The minimum size of the void phase is indicated in blue, while the minimum size of the solid phase in magenta. The table also reports the 3 parameters required to impose the desired minimum length scales, i.e.~$\Threshold{ero}$, $\Threshold{dil}$, and $r_\mathrm{fil}$. These parameters have been obtained using the analytical method implemented in the MATLAB code provided with this paper (\texttt{SizeSolution.m}). It is recalled that all the examples assume $\Threshold{int}=0.5$.

%%%%%%%%%%%%%%%%%%%%%%%%%%%%%%%%%%%
%%%%%%%%%%%%%%%%%%%%%%%%%%%%%%%%%%%

%====================================
%====================================
\begin{table*}
\centering
    \begin{tabular}{c c c c }
      \toprule
      \multirow{3}{*}{\large Mesh}&  \multicolumn{3}{c}{ \large $\rminsolid{int}/\rminvoid{int}$} 
      \vspace{1mm} \\
      \cline{2-4}
      &  
      \large $1/2$ 
      \vspace{0.5mm}  
      & \large $1/1$ 
      & \large $2/1$ \\
      & using [$\Threshold{ero}$, $\Threshold{dil}$] = [0.65, 0.05]
      & using [$\Threshold{ero}$, $\Threshold{dil}$] = [0.70, 0.30]
      & using [$\Threshold{ero}$, $\Threshold{dil}$] = [0.80, 0.42]
      \vspace{1mm}\\
      \hline
      \multirow{2}{*}[7em]{$100\times100$}  & 
      \includegraphics[width=0.25\linewidth]{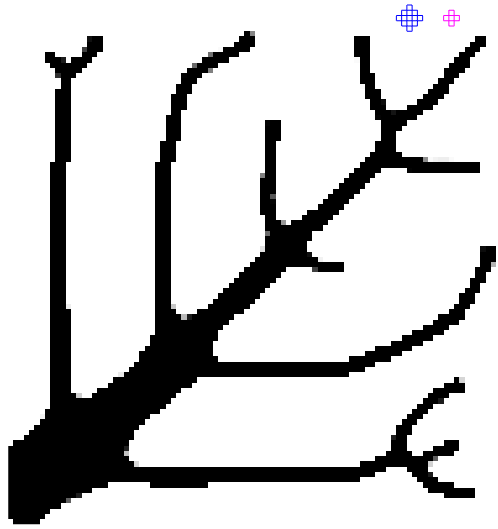}
        \label{fig:Thermal100_2/3}& 
      \includegraphics[width=0.25\linewidth]{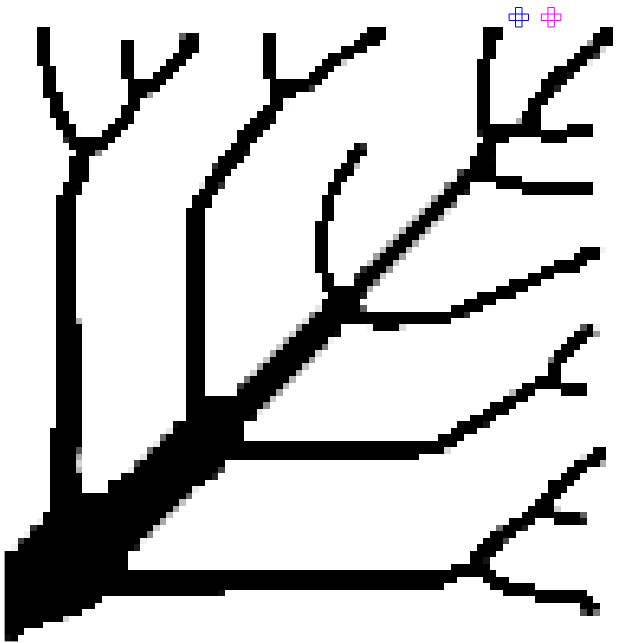}
        \label{fig:Thermal100_2/2}&
      \includegraphics[width=0.25\linewidth]{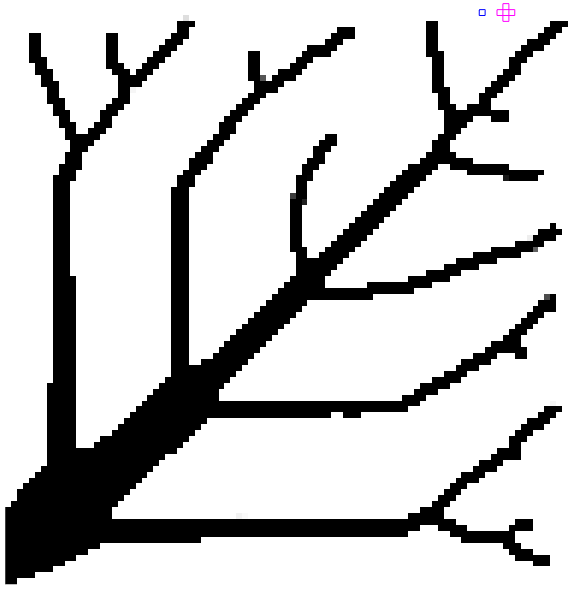}
        \label{fig:Thermal100_3/2}
      \\
      	&\vspace{0.5mm}
		$r_\mathrm{min.Solid}^\mathrm{int} = 1$, $r_\mathrm{fil}= 2.58$.	
		\vspace{0.7mm}
		&
		$r_\mathrm{min.Solid}^\mathrm{int} = 1$, $r_\mathrm{fil}= 2.24$.    
		&
		$r_\mathrm{min.Solid}^\mathrm{int} = 1$, $r_\mathrm{fil}= 1.81$.
      \\
      \hline
       \multirow{2}{*}[7em]{$200\times200$} & 
      \includegraphics[width=0.25\linewidth]{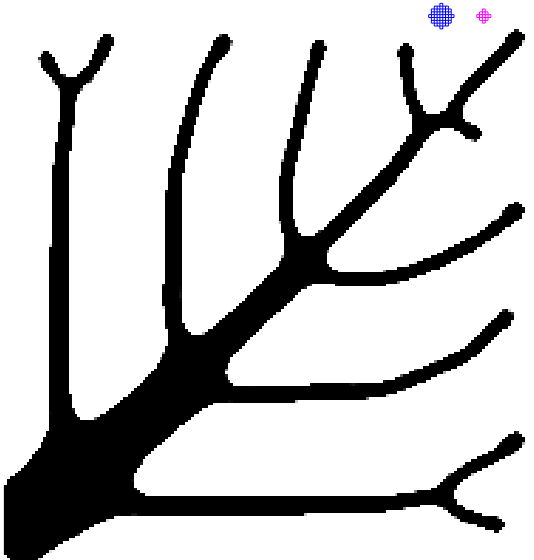}
        \label{fig:Thermal200_2/3}&
      \includegraphics[width=0.25\linewidth]{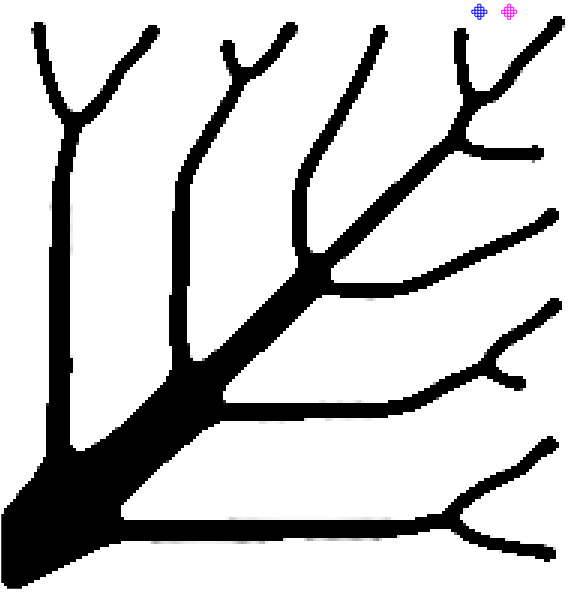}
		\label{fig:Thermal200_2/2} &
	  \includegraphics[width=0.25\linewidth]{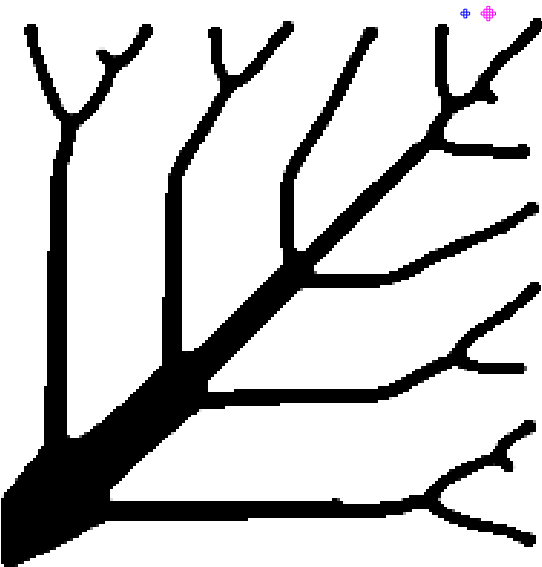}
	    \label{fig:Thermal200_3/2}
	    \\
		&\vspace{0.5mm}
		$r_\mathrm{min.Solid}^\mathrm{int} = 2$, $r_\mathrm{fil}= 5.16$.
		\vspace{0.7mm}
		&
		$r_\mathrm{min.Solid}^\mathrm{int} = 2$, $r_\mathrm{fil}= 4.47$.    
		&
		$r_\mathrm{min.Solid}^\mathrm{int} = 2$, $r_\mathrm{fil}= 3.62$.    
	    \\
      \hline
       \multirow{2}{*}[7em]{$400\times400$} &
      \includegraphics[width=0.25\linewidth]{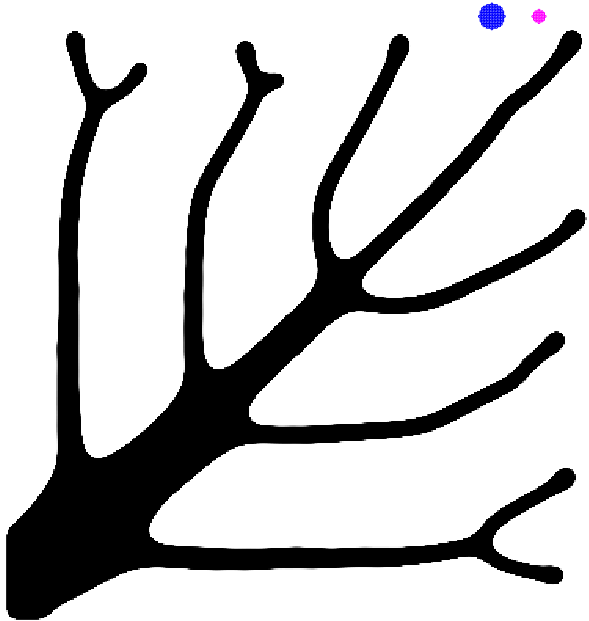}
        \label{fig:Thermal400_2/3}& \includegraphics[width=0.25\linewidth]{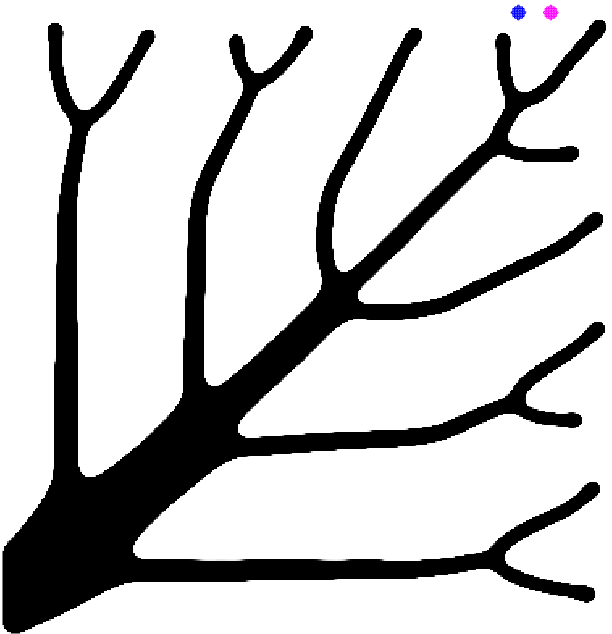}
        \label{fig:Thermal400_2/2}& 
      \includegraphics[width=0.25\linewidth]{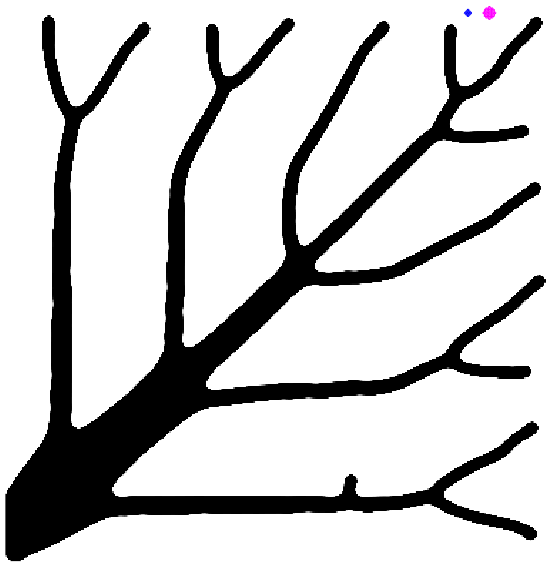}
        \label{fig:Thermal400_3/2}
        \\
        &\vspace{0.5mm}
		$r_\mathrm{min.Solid}^\mathrm{int} = 4$, $r_\mathrm{fil}= 10.32$.	
		\vspace{0.7mm}
		&
		$r_\mathrm{min.Solid}^\mathrm{int} = 4$, $r_\mathrm{fil}= 9.94$.    
		&
		$r_\mathrm{min.Solid}^\mathrm{int} = 4$, $r_\mathrm{fil}= 7.24$.
        \\
      \bottomrule
    \end{tabular}
\caption{{Optimized designs for the heat exchanger problem. The imposed volume constraint is $V^*_\mathrm{min}=20\%$.}}
\label{Tab:9_HeatExchangers}
\end{table*}
%The nine solutions are obtained from the combination of three different discretizations and three different minimum length scales.
%====================================
%====================================

\begin{figure}[t!]
		\captionsetup[subfigure]{labelformat=empty}
    	\centering
     	\begin{subfigure}{0.98\linewidth}
		\centering \includegraphics[width=0.98\linewidth]{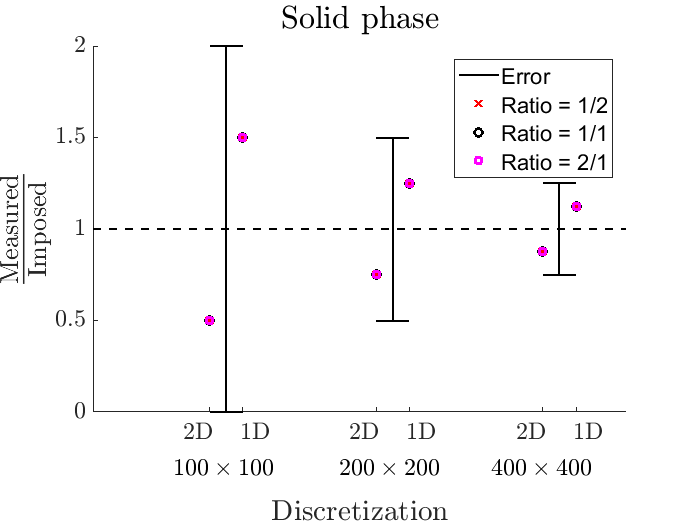}
		\end{subfigure}		
		\vspace{4mm}\\
		\begin{subfigure}{0.95\linewidth}
		\centering \includegraphics[width=0.95\linewidth]{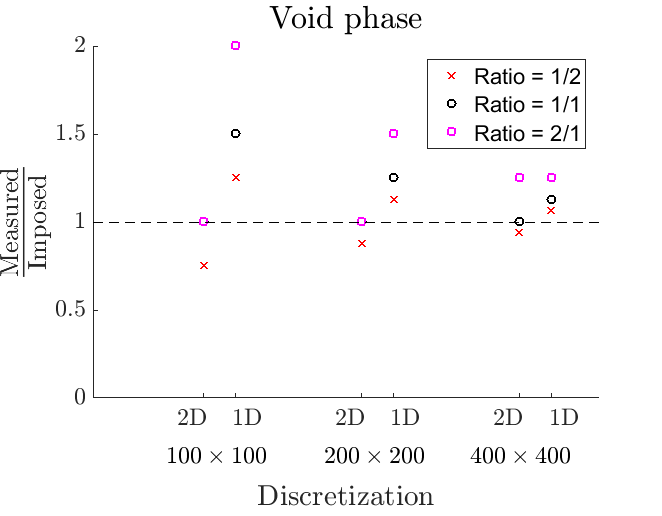}
		\end{subfigure}	
		\\
		\begin{minipage}{0.01\textwidth}
			\subcaption{\label{Fig:DiscreteInfluence_a}}
		\end{minipage}
		~
		\begin{minipage}{0.01\textwidth}
			\subcaption{\label{Fig:DiscreteInfluence_b}}
		\end{minipage}
		\vspace{-126mm}\\
		\hspace{-78mm}\footnotesize{(a)} 
		\vspace{62mm}\\
		\hspace{-78mm}\footnotesize{(b)} 
		\vspace{54mm}\\
		\caption{The graphs summarize the measured minimum size of (a) the solid phase and (b) the void phase from the designs in Table \ref{Tab:9_HeatExchangers}. For each discretization, two set of results are provided, labeled as $2D$ and $1D$. The $2D$ represents the graphical measurement normalized with respect to the imposed value. The $1D$ represents the minimum size obtained from the numerical method normalized with respect to the imposed value.}
		\label{Fig:DiscreteInfluence}
\end{figure}

The measured minimum size of the solid and void phases are reported on Figs.~\ref{Fig:DiscreteInfluence_a} and \ref{Fig:DiscreteInfluence_b}, respectively. To construct these graphs, the measured minimum sizes are normalized with respect to the desired minimum sizes, and the values are placed in the ordinate coordinate. This is done for each discretization (which determines the abscissa) and for each length scale. This procedure is also carried out with results obtained from the numerical method \citep{Wang2011}, which is executed using representative discretizations for each case. The values obtained from the numerical method are labeled as \textit{1D} in the graphs of Fig.~\ref{Fig:DiscreteInfluence}, while the values measured from the optimized designs are labeled as \textit{2D}. The error bars in Fig.~\ref{Fig:DiscreteInfluence_a} illustrate the rounding error present in the analytical method (see Section \ref{Sec:rounding_error} for definition). Rounding error lines have not been plotted in Fig.~\ref{Fig:DiscreteInfluence_b} because the length scale ratios have different rounding errors.

As a general observation, we can point out from Fig.~\ref{Fig:DiscreteInfluence} that mesh refinement reduces the error between the measured minimum size and the desired minimum size (which is imposed through the set of parameters provided by the analytical method). This is consistent with the observations made in the previous section, where a continuous and uni-dimensional design domain was used. In addition, the predicted error for the analytical method (the error bars) are also consistent with the measured sizes, which validates the scope and limitations of the analytical equations discussed in Section \ref{sec:5}.

Regarding the graph in Fig.~\ref{Fig:DiscreteInfluence_a}, it can be mentioned that in the nine 2D designs, the error between the measured and desired minimum sizes is half a finite element. This error is relatively large in the coarse discretization (50$\%$ error), and small in the fine discretization (12.5$\%$ error). However, despite the fact that the analytical method is not exact, it seems to be accurate enough in the examples examined, since in all cases the error is the same, half a finite element. It is interesting to note that the numerical method that assumes a discrete 1D domain, estimates a minimum size that differs by half a finite element with respect to the imposed value, which matches the measured error. However, the numerical method provides a minimum size of the solid phase larger than the measured one. This is probably due to the fact that the condition of robustness refers to one finite element and not to an infinitesimal size, as discussed in Section \ref{Sec:Infinitesimal_Size}.

On the other hand, the graph built for the void phase (Fig.~\ref{Fig:DiscreteInfluence_b}) shows a different pattern. The measured minimum size of the void phase is often equal to the desired one and even bigger for some designs. This could be explained by the simple fact that the chosen test case does not present a geometric singularity over the minimum size of the void phase, therefore it is not possible to guarantee that the smallest re-entrant corner of the design will indeed correspond to the minimum size imposed by the robust formulation. The representative case might be that where $r_\mathrm{min.Void}^\mathrm{int}$ is chosen twice the size $r_\mathrm{min.Solid}^\mathrm{int}$ (ratio 1/2). In such a case, the measured error corresponds to half a finite element smaller than the desired value, that is, the same result as for the solid phase.

\input{ErosionAndDilation}

\subsection{Remark on the simplified robust formulation}

Before concluding this manuscript we devote a discussion concerning the robust formulation. It is widely known that the robust formulation of the optimization problem can be simplified when the objective function is monotonically dependent on the volume fraction. The most widespread case in the literature is the compliance minimization problem, where the intermediate and dilated fields have less compliance than the eroded design and therefore can be removed from the objective function without compromising the robustness of the formulation. In this case, the objective function is evaluated only for the eroded design, with the intermediate and dilated fields remaining solely for formulating design constraints. When the volume restriction is the only constraint included in the optimization problem, the constraint can be formulated using the intermediate design, assuming that this design field is the one intended for manufacturing. However, it has been mentioned in the literature that evaluating the volume restriction in the dilated design promotes convergence to better optimums since numerical instabilities are prevented.
In this section we add another reason that has not been mentioned so far (to the best of the authors' knowledge). If the dilated design is not included in the volume restriction, then it is not possible to ensure the control over the minimum size of the void phase. To explain this statement, consider the following two robust topology optimization formulations for the thermal compliance minimization problem:  

\newpage
\begin{flalign*}
\hspace{14mm} \text{\small{P.I}} \hspace{33mm} \text{\small{P.II}} 
\end{flalign*}
\vspace{-8mm}
$$
\hspace{-10mm} \overbrace{\hspace{30mm}} 
\hspace{4mm} \overbrace{\hspace{37mm}} 
$$
\vspace{-8mm}
\begin{align} \label{EQ:OPTIS_Non_Robust}
  		{\min_{\bm{\rho}}} & \quad c({\ProjField{ero}}) 
  		&  		
  		{\min_{\bm{\rho}}} & \quad c({\ProjField{ero}})  
  		\nonumber \\
  		\mathrm{s.t. :} & \quad \mathbf{v}^\intercal {\ProjField{int}} \leq V^*_\mathrm{int} \:,
  		&
  		\mathrm{s.t. :} & \quad \mathbf{v}^\intercal {\ProjField{dil}} \leq V^*_\mathrm{dil}(V^*_\mathrm{int})
  		\\
  		& \quad	 0 \leq {\rho_i} \leq1
  		&
  		& \quad 0 \leq {\rho_i} \leq1
  		\nonumber
\end{align}

The two optimization problems, P.I and P.II in \linebreak Eq.~\eqref{EQ:OPTIS_Non_Robust}, are formulated under the robust design approach, but P.I evaluates the volume constraint directly in the intermediate design while P.II does it through the dilated design, whose upper bound $V^*_\mathrm{dil}$ is scaled according to $V^*_\mathrm{int}$. From Eq.~\eqref{EQ:OPTIS_Non_Robust} it is clear that P.I$({\ProjField{ero}},{\ProjField{int}})$ and P.II$({\ProjField{ero}},{\ProjField{int}},{\ProjField{dil}})$, or alternatively, P.I$({\Threshold{ero}},{\Threshold{int}})$ and P.II$({\Threshold{ero}},{\Threshold{int}},{\Threshold{dil}})$. Therefore, P.I is not influenced by the dilation threshold nor the dilated design. 

\begin{figure}
	\captionsetup{width=1.00\linewidth}
    \centering	
    \begin{subfigure}[t]{0.48\linewidth}
    	\centering
		\includegraphics[width=0.98\linewidth]{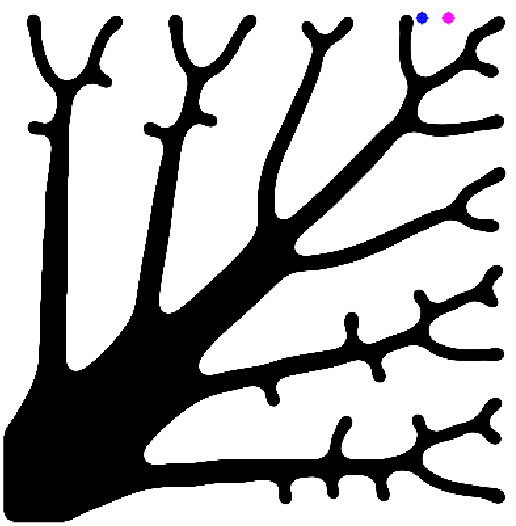}
		\caption{P.II, ${\Threshold{dil}} = 0.40$.}
		\label{fig:dil40}
	\end{subfigure}
    ~
	\begin{subfigure}[t]{0.48\linewidth}
    	\centering
    	\includegraphics[width=0.98\linewidth]{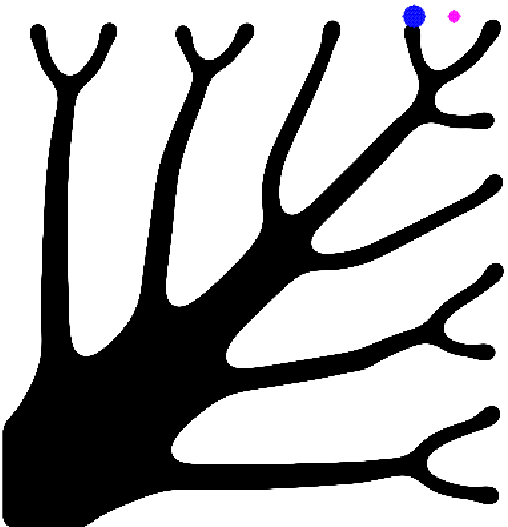}
		\caption{P.II, ${\Threshold{dil}} = 0.14$.}
		\label{fig:dil14}
	\end{subfigure}
	\vspace{3mm}\\
	\begin{subfigure}[t]{0.48\linewidth}
    	\centering
		\includegraphics[width=0.98\linewidth]{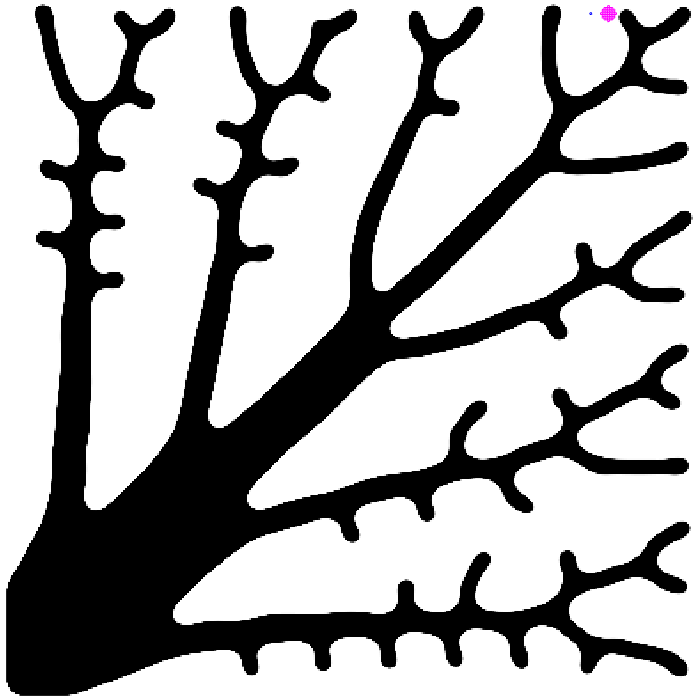}
		\caption{P.II, ${\Threshold{dil}} = 0.50$.}
		\label{fig:dil50}
		~%\vspace{3mm}\\
	\end{subfigure}
	\begin{subfigure}[t]{0.48\linewidth}
    	\centering
		\includegraphics[width=0.98\linewidth]{Figures/400Sol4Void0.png}
		\caption{P.I for any ${\Threshold{dil}}$.}
		\label{fig:dil025Bis}
	\end{subfigure}

	\caption{Heat exchanger design problem using two different variations of the robust formulation. P.I evaluates the volume constraint directly in the intermediate design, while P.II does it through the dilated design. Here, ${\Threshold{ero}} = 0.60$ and $r_\mathrm{min.Solid}^\mathrm{int}=4$ elements.}
	\label{fig:HeatSinkDil}
\end{figure}

We recall that the condition of robustness imposed for the void phase involves the dilated design, which has to project a cavity with at least one void element to be present in the 3 fields that constitute the robust formulation. The influence of the dilated design on the minimum size of the void phase can be seen graphically in Fig.~\ref{fig:fig_1_a} by considering a fixed erosion threshold and different dilation thresholds. For example, for the set $[{\Threshold{dil}}$, ${\Threshold{ero}}]$ consider the points $[0.14 \:,\: 0.60]$ and $[0.40 \;,\; 0.60]$. The first point corresponds to a length scale where the minimum size of the void phase is equal to that of the solid phase, while the second point imposes the size of the void phase twice that of the solid phase. In the following, problems P.I and P.II are solved using the two sets of thresholds indicated previously. To help to see the differences, the volume fraction ratio has been increased to $0.3$. These sets of parameters are resumed in Table \ref{tab:SizeSolutionVoid}. The results are summarized in Fig.~\ref{fig:HeatSinkDil}.

Clearly, the results obtained using P.II show a length scale consistent with the expected one, as the design in Fig.~\ref{fig:dil40} features bigger reentrant corners than the design in Fig.~\ref{fig:dil14}. However, for any value of ${\Threshold{dil}}$, the result from P.I is always the same and corresponds to that shown in Fig.~\ref{fig:dil025Bis}. The interpretation that can be made of the P.I problem is that it uses a dilation threshold equal to the intermediate one, i.e.~${\Threshold{dil}}= {\Threshold{int}} = 0.5$. This can be corroborated by solving P.II with the set $[{\Threshold{dil}}\:,\: {\Threshold{ero}}]$ equal to $[0.50, 0.60]$, whose result is shown in Fig.~\ref{fig:dil50}.

\begin{table}
    \centering
    \begin{tabular}{c c c c c}
    \toprule
      $\rminsolid{int}$ & $\rminvoid{int}$ &  $\rfil$ &  $\eta_e$ &  $\eta_d$
      \vspace{1mm}\\
      \hline
      $4$ & $0$& $12.65$ & $0.60$& $0.50$\\

      $4$ & $4$& $12.65$ & $0.60$& $0.40$\\

      $4$ & $8$& $12.65$ & $0.60$ & $0.14$ \\
      \hline
    \end{tabular}
    \caption{Parameters to compare the robust formulation and its simplified form.}
    \label{tab:SizeSolutionVoid}
\end{table}

\section{Conclusion}\label{sec:7}

The robust topology optimization formulation based on uniform manufacturing errors has gained increasing acceptance in the topology optimization community. This is mainly due to its ability to control the minimum size of both the solid and void phases, and its potential to be combined with other topology optimization approaches. Despite the increasing popularity of the formulation, no method was yet available to easily obtain the filter and projection parameters that produce the desired minimum length scales. This need encouraged us to further develop the analytical method proposed by \citet{Qian2013}. The scope and limitations of this method were assessed using the numerical method of \citet{Wang2011} and a set of 2D design results from two topology optimization problems, the thermal compliance minimization problem and the non-linear force inverter.

In addition to providing a fast and effective way to obtain the parameters that impose the desired minimum length scale, this work shows that to obtain simultaneous control over the minimum sizes of the solid and void phases, it is necessary to involve the 3 fields (eroded, intermediate and dilated) in the robust topology optimization problem. For example, for the compliance minimization problem subject to a volume restriction, it is known that intermediate and dilated designs can be excluded from the objective function, but the volume restriction has to be applied to the dilated design in order to involve all 3 designs in the formulation.

%=============================================================================================
\section{Replication of results}\label{sec:8}
This manuscript contains two MATLAB codes as supplementary material. The first is called \texttt{SizeSolution.m} and provides a list of filter and projection parameters that impose user defined minimum length scales. The second is called \texttt{NumericalSolution.m} and builds the graphs in Fig.~\ref{fig:fig_1} using the numerical method proposed by \citet{Wang2011}.

\begin{acknowledgements}
The authors acknowledge the research project FAFIL (\textit{Fabrication Additive par D\'{e}p\^{o}t de Fil}), funded by INTERREG and the European Regional Development Fund (ERDF).
\end{acknowledgements}

\section*{Conflict of interest}
On behalf of all authors, the corresponding author states that there is no conflict of interest.

\bibliographystyle{spbasic}      % basic style, author-year citations
\bibliography{References}   % name your BibTeX data base
%=============================================================================================

\end{document}

%% file: ErosionAndDilation.tex
\subsection{Erosion and dilation distances}

In the following example we illustrate  the use of the erosion and dilation distances that are provided for a desired minimum length scale. To this end, we use maximum size constraints, where the erosion and dilation distances are essential information to impose a consistent length scale in the robust formulation \citep{Fernandez2020}. For the sake of completeness of the manu\-script, the non-linear force inverter is considered in this illustrative example. The topology optimization problem of the  non-linear  force  inverter  including  maximum  size constraint reads as follow :

\begin{equation}
\begin{split}
    \text{min} \quad & \text{max}\left( \mathrm{c}({\ProjField{ero}}),\; \mathrm{c}({\ProjField{int}})\;, \mathrm{c}({\ProjField{dil}}) \right)   \\
    \text{s.t.:} \quad & \mathbf{v}^\intercal \: {\ProjField{dil}} \leq V^*_\mathrm{dil}(V^*_\mathrm{int}) \\
                & \mathrm{G_{ms}}(\ProjField{dil}) \leq 0 \\
                & 0 \leq \rho_i \leq 1 \;\; , \;\; i=1,...,N\;\;,
\end{split}
\label{Eq:NLInverter}
\end{equation}

\noindent where $\mathrm{G_{ms}}$ is the maximum size constraint defined exactly as in \citep{Fernandez2020}. Therefore, $\mathrm{G_{ms}}$ represents a \textit{p-mean} aggregation function that gathers local maximum size restrictions. As in the cited work, the $p$ aggregation exponent is set at 100. 

Regarding the optimization parameters, these represent the implementation of \cite{Wang2014}, i.e. the Heaviside parameter $\beta$ is initialized at 1 and is increased by 1 every 20 iterations until a value $\beta = 16$ is reached. Then, 20 more iterations are carried out with a $\beta = 32$. The SIMP penalty parameter is set at 3. We have found that this setting of parameters works well for introducing maximum size restrictions into the non-linear force inverter formulated under the robust design approach. 

To impose minimum and maximum length  scales, the maximum size constraint $\mathrm{G_{ms}}$ is applied on the dilated design \citep{Fernandez2020}. To do so, the regions where the local maximum size constraints are applied must be scaled according to the dilation distance. For instance, if the desired maximum size in the intermediate design is defined by a circle of radius $r_\mathrm{max.Solid}^\mathrm{int}$, the maximum size constraint should be formulated for the dilated design using a circle of radius:
\begin{equation}
    r_\mathrm{max.Solid}^\mathrm{dil} = r_\mathrm{max.Solid}^\mathrm{int} + t_\mathrm{dil}
    \label{Eq:rMaxDil}
\end{equation}

Equation \eqref{Eq:rMaxDil} explains the need of knowing the dilation distance when imposing maximum size restrictions in the robust formulation. As mentioned previously, this information can be easily obtained from the graphs generated by the analytical equations. For example, in the following, the half force inverter depicted in Fig. \ref{fig:DesignDomains} is solved. The design domain is discretized into $200 \times 100$ quadrilateral finite elements. The minimum size of the solid phase is set as $r_\mathrm{min.Solid}^\mathrm{int} = 2$ elements, while the maximum size is set as $r_\mathrm{max.Solid}^\mathrm{int} = 3$  elements. Two different values for the minimum size of the void phase are  chosen, $r_\mathrm{min.Void}^\mathrm{int} = 2$ and $r_\mathrm{min.Void}^\mathrm{int} = 3$ elements. The filter and projection parameters used to impose the desired length scale are reported in Table \ref{tab:SizeSolutionBis}. In all cases, the volume constraint is set to $25\%$.
\begin{table}[b]
    \centering
    \begin{tabular}{c c c c c c}
    \toprule
       $r_\mathrm{min.Void}^\mathrm{int}$ & $\rfil$ &  $\eta_\mathrm{ero}$ &  $\eta_\mathrm{dil}$ &  $t_\mathrm{dil}$ & $t_\mathrm{ero}$
       \vspace{1mm}\\
      \hline
       $2$ & $4.47$ & $0.70$& $0.30$ & $1.03$ & $1.03$\\

      $3$ &$4.47$ & $0.70$& $0.11$ & $2.41$ & $1.03$  \\
      \hline
    \end{tabular}
    \caption{Parameters to impose $r_\mathrm{min.Solid}^\mathrm{int} = 2$ on the design with $\beta = 32$.}
    \label{tab:SizeSolutionBis}
\end{table}

To obtain the dilation distance $t_\mathrm{dil}$, the projection thresholds that define the minimum length scale must be defined. These values can be graphically obtained from Fig.~\ref{fig:fig_1}, or from the attached MATLAB code \linebreak (\texttt{SizeSolution.m}). For example, using this code, the sets of parameters in Table \ref{tab:SizeSolutionBis} are obtained. The dilation distance is also reported there. For the minimum sizes $r_\mathrm{min.Void}^\mathrm{int} = 2$ and $r_\mathrm{min.Void}^\mathrm{int} = 3$ elements, the dilation distances are $t_\mathrm{dil}=1$ and $t_\mathrm{dil}=2$ elements, respectively (the numbers have been rounded to the nearest integer). 

\begin{figure}
	\captionsetup{width=1\linewidth}
    \centering
		\includegraphics[width=0.75\linewidth , angle =0]{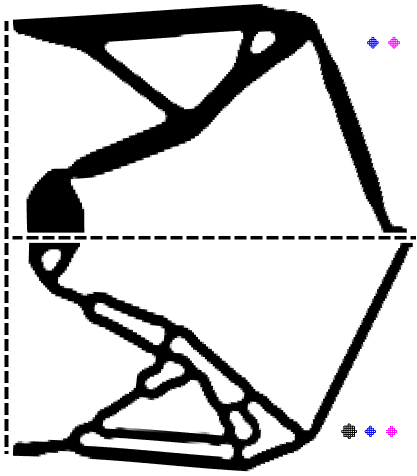}
\caption{Nonlinear force inverter without (upper half) and with (lower half) maximum size constraints. The minimum length scales are $r_\mathrm{min.Solid}^\mathrm{int}$ = $r_\mathrm{min.Void}^\mathrm{int} $ = 2 elements. The maximum size is $r_\mathrm{max.Solid}^\mathrm{int}=3$ elements.}
	\label{fig:ForceInvertersMaxVoid2}
\end{figure}

\begin{figure}
	\captionsetup{width=1.00\linewidth}
    \centering	
		\includegraphics[width=0.75\linewidth, angle =0]{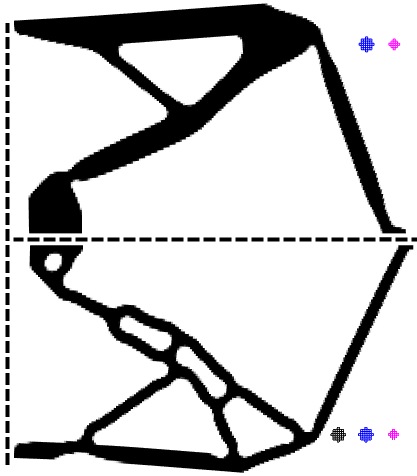}
	\caption{Nonlinear force inverter without (upper half) and with (lower half) maximum size constraints. The minimum length scales are $r_\mathrm{min.Solid}^\mathrm{int}$ = 2 elements and $r_\mathrm{min.Void}^\mathrm{int}$ = 3 elements. The prescribed maximum size is $r_\mathrm{max.Solid}^\mathrm{int}=3$ elements.}
	\label{fig:Inverter23max3}
\end{figure}

The results are shown in Figs.~\ref{fig:ForceInvertersMaxVoid2} and \ref{fig:Inverter23max3}. Each figure reports 2 optimized designs, which are obtained with and without maximum size restrictions. To facilitate the visual comparison between results, the designs are placed with respect to the symmetry axis and are reported in their deformed configuration. The imposed length scale is also reported graphically next to each design. As in the previous examples, the blue and magenta circles represent the minimum size of the void and solid phase, respectively. The black circle represents the maximum size desired for the intermediate design (which is imposed through the dilated design)\footnote{The maximum size is actually imposed using an annular region \citep{Fernandez2020}, but for illustrative purposes, a circular region is drawn.}. 

The results are consistent with the imposed length scales. In the 4 designs reported in Figs.~\ref{fig:ForceInvertersMaxVoid2} and \ref{fig:Inverter23max3}, the minimum size of the solid and void phases are met with a half-finite element of error (which agrees with the analytical rounding error). Regarding the maximum size, this also matches the imposed value. However, due to the inherent drawbacks of the aggregation function and the strong non-linearity of the force inverter design problem, there are local regions where the imposed maximum size restriction is not met, such as the horizontal bar in the design of Fig.~\ref{fig:Inverter23max3}.    

This force inverter test case shows the usefulness of the proposed equations and provided codes, since they allow to quickly obtain the filter and projection parameters, and the dilation distance that allow to impose the desired length scales. 

%% file: Article.bbl
\begin{thebibliography}{27}
\providecommand{\natexlab}[1]{#1}
\providecommand{\url}[1]{{#1}}
\providecommand{\urlprefix}{URL }
\expandafter\ifx\csname urlstyle\endcsname\relax
  \providecommand{\doi}[1]{DOI~\discretionary{}{}{}#1}\else
  \providecommand{\doi}{DOI~\discretionary{}{}{}\begingroup
  \urlstyle{rm}\Url}\fi
\providecommand{\eprint}[2][]{\url{#2}}

\bibitem[{Andreasen et~al(2020)Andreasen, Elingaard, and Aage}]{Andreasen2020}
Andreasen CS, Elingaard MO, Aage N (2020) Level set topology and shape
  optimization by density methods using cut elements with length scale control.
  Structural and Multidisciplinary Optimization pp 1--23

\bibitem[{Bends{\o}e and Kikuchi(1988)}]{BensoeKikuchi1988}
Bends{\o}e MP, Kikuchi N (1988) Generating optimal topologies in structural
  design using a homogenization method. Computer methods in applied mechanics
  and engineering 71(2):197--224

\bibitem[{Bendsøe(1989)}]{Bendsoe1989}
Bendsøe M (1989) Bendsoe, m.p.: Optimal shape design as a material
  distribution problem. structural optimization 1, 193-202. Structural
  Optimization 1:193--202

\bibitem[{Bourdin(2001)}]{Bourdin2001}
Bourdin B (2001) Filters in topology optimization. International journal for
  numerical methods in engineering 50(9):2143--2158

\bibitem[{Bruns and Tortorelli(2001)}]{Bruns2001}
Bruns TE, Tortorelli DA (2001) Topology optimization of non-linear elastic
  structures and compliant mechanisms. Computer methods in applied mechanics
  and engineering 190(26-27):3443--3459

\bibitem[{Chen and Chen(2011)}]{Chen2011}
Chen S, Chen W (2011) A new level-set based approach to shape and topology
  optimization under geometric uncertainty. Structural and Multidisciplinary
  Optimization 44(1):1--18

\bibitem[{Christiansen et~al(2015)Christiansen, Lazarov, Jensen, and
  Sigmund}]{Christiansen2015}
Christiansen R, Lazarov B, Jensen J, Sigmund O (2015) Creating geometrically
  robust designs for highly sensitive problems using topology optimization:
  Acoustic cavity design. Structural and Multidisciplinary Optimization
  52:737--754

\bibitem[{Clausen and Andreassen(2017)}]{Clausen2017}
Clausen A, Andreassen E (2017) On filter boundary conditions in topology
  optimization. Structural and Multidisciplinary Optimization 56(5):1147--1155

\bibitem[{{da Silva} et~al(2019){da Silva}, Beck, and Sigmund}]{DASILVA2019_2}
{da Silva} GA, Beck AT, Sigmund O (2019) Topology optimization of compliant
  mechanisms with stress constraints and manufacturing error robustness.
  Computer Methods in Applied Mechanics and Engineering 354:397 -- 421

\bibitem[{Fern{\'a}ndez et~al(2020)Fern{\'a}ndez, Yang, Koppen, Alarc{\'o}n,
  Bauduin, and Duysinx}]{Fernandez2020}
Fern{\'a}ndez E, Yang Kk, Koppen S, Alarc{\'o}n P, Bauduin S, Duysinx P (2020)
  Imposing minimum and maximum member size, minimum cavity size, and minimum
  separation distance between solid members in topology optimization. Computer
  Methods in Applied Mechanics and Engineering 368:113,157

\bibitem[{Fern{\'a}ndez et~al(2021)Fern{\'a}ndez, Ayas, Langelaar, and
  Duysinx}]{Fernandez2021}
Fern{\'a}ndez E, Ayas C, Langelaar M, Duysinx P (2021) Topology optimization
  for large-scale additive manufacturing: Generating designs tailored to the
  deposition nozzle size (Under Review.)

\bibitem[{Kumar and Fern{\'a}ndez(2021)}]{Kumar2021}
Kumar P, Fern{\'a}ndez E (2021) A numerical scheme for filter boundary
  conditions in topology optimization on regular and irregular meshes (Under
  Review, arXiv:2101.01122v1)

\bibitem[{Lazarov and Sigmund(2011)}]{Lazarov2011}
Lazarov BS, Sigmund O (2011) Filters in topology optimization based on
  helmholtz-type differential equations. International Journal for Numerical
  Methods in Engineering 86(6):765--781

\bibitem[{Pedersen and Allinger(2006)}]{Pedersen2006}
Pedersen C, Allinger P (2006) Industrial Implementation and Applications of
  Topology Optimization and Future Needs, vol 137, Springer, pp 229--238

\bibitem[{Pellens et~al(2018)Pellens, Lombaert, Lazarov, and
  Schevenels}]{Pellens2018}
Pellens J, Lombaert G, Lazarov B, Schevenels M (2018) Combined length scale and
  overhang angle control in minimum compliance topology optimization for
  additive manufacturing. Structural and Multidisciplinary Optimization

\bibitem[{Qian and Sigmund(2013)}]{Qian2013}
Qian X, Sigmund O (2013) Topological design of electromechanical actuators with
  robustness toward over-and under-etching. Computer Methods in Applied
  Mechanics and Engineering 253:237--251

\bibitem[{Sigmund(1997)}]{Sigmund1997}
Sigmund O (1997) On the design of compliant mechanisms using topology
  optimization. Journal of Structural Mechanics 25(4):493--524

\bibitem[{Sigmund(2009)}]{Sigmund2009}
Sigmund O (2009) Manufacturing tolerant topology optimization. Acta Mechanica
  Sinica 25(2):227--239

\bibitem[{Sigmund and Maute(2013)}]{Sigmund2013}
Sigmund O, Maute K (2013) Topology optimization approaches. Structural and
  Multidisciplinary Optimization 48(6):1031--1055

\bibitem[{Silva et~al(2020)Silva, Beck, and Sigmund}]{daSilva2020}
Silva G, Beck A, Sigmund O (2020) Topology optimization of compliant mechanisms
  considering stress constraints, manufacturing uncertainty and geometric
  nonlinearity. Computer Methods in Applied Mechanics and Engineering
  365:112,972

\bibitem[{Wang et~al(2011)Wang, Lazarov, and Sigmund}]{Wang2011}
Wang F, Lazarov BS, Sigmund O (2011) On projection methods, convergence and
  robust formulations in topology optimization. Structural and
  Multidisciplinary Optimization 43(6):767--784

\bibitem[{Wang et~al(2011b)Wang, Jensen, and Sigmund}]{Wang2011b}
Wang F, Jensen J, Sigmund O (2011b) Robust topology optimization of photonic
  crystal waveguides with tailored dispersion properties. JOSA B 28:387--397

\bibitem[{Wang et~al(2014)Wang, Lazarov, Sigmund, and Jensen}]{Wang2014}
Wang F, Lazarov BS, Sigmund O, Jensen JS (2014) Interpolation scheme for
  fictitious domain techniques and topology optimization of finite strain
  elastic problems. Computer Methods in Applied Mechanics and Engineering
  276:453 -- 472

\bibitem[{Xu et~al(2010)Xu, Cai, and Cheng}]{Xu2010}
Xu S, Cai Y, Cheng G (2010) Volume preserving nonlinear density filter based on
  heaviside functions. Structural and Multidisciplinary Optimization
  41(4):495--505

\bibitem[{Yan et~al(2018)Yan, Wang, and Sigmund}]{Yan2018}
Yan S, Wang F, Sigmund O (2018) On the non-optimality of tree structures for
  heat conduction. International Journal of Heat and Mass Transfer 122:660--680

\bibitem[{Zhou et~al(2011)Zhou, Fleury, Patten, Stannard, Mylett, and
  Gardner}]{Zhou2011}
Zhou M, Fleury R, Patten S, Stannard N, Mylett D, Gardner S (2011) Topology
  optimization-practical aspects for industrial applications. In: 9th World
  Congress on Structural and Multidisciplinary Optimization

\bibitem[{Zhu et~al(2016)Zhu, Zhang, and Xia}]{Zhu2016}
Zhu JH, Zhang WH, Xia L (2016) Topology optimization in aircraft and aerospace
  structures design. Archives of Computational Methods in Engineering
  23:595--622

\end{thebibliography}
